\documentclass{article}
\usepackage{lineno}
\usepackage[a4paper,margin=1in,top=3cm, bottom=3cm]{geometry}
\usepackage[a4paper]{geometry}
\usepackage[utf8]{inputenc}
\usepackage{amsmath,amsfonts,amssymb,amsthm,bbm,mathrsfs}
\usepackage[shortlabels]{enumitem}
\usepackage{xcolor}
\usepackage{graphicx}
\usepackage{appendix}
\usepackage{setspace}
\usepackage{fancyhdr}
\usepackage[
    backref=page
]{hyperref}

\hypersetup{
    colorlinks=true,
    linkcolor=blue,
    filecolor=magenta,      
    urlcolor=blue,
    citecolor=blue,
    pdfpagemode=FullScreen,
}

\newtheorem{thm}{Theorem}[section]

\newtheorem{rmk}{Remark}[section]

\newtheorem{prop}{Proposition}[section]
\newtheorem{lm}{Lemma}[section]

\numberwithin{equation}{section}

\makeatletter
\def\@setauthors{
  \begingroup
  \def\thanks{\protect\thanks@warning}%
  \trivlist
  \centering\footnotesize \@topsep30\p@\relax
  \advance\@topsep by -\baselineskip
  \item\relax
  \author@andify\authors
  \def\\{\protect\linebreak}%
  \authors%
  \ifx\@empty\contribs
  \else
    ,\penalty-3 \space \@setcontribs
    \@closetoccontribs
  \fi
  \endtrivlist
  \endgroup
}
\makeatother
\fancypagestyle{firstpage}{
\fancyhf{} % clear headers and footers

\fancyfoot[L]{%
\footnotesize
\rule{3cm}{0.3pt}\\ 
This article is based upon work from the project PRIN2022 D53D23005580006 ``Elliptic and parabolic problems, heat kernel estimates and spectral theory''. The second and third authors are members of the ``Gruppo Nazionale per l’Analisi Matematica, la Probabilità e le loro Applicazioni (GNAMPA)'' of the Istituto Nazionale di Alta Matematica (INdAM).
  }
}

\title{\textbf{Controllability of Forward Stochastic Reaction--Convection--Diffusion Systems with Cascade Structure}}

\author{Abdellatif Elgrou$^\dag\,^*$,\,\,Federica Gregorio$^\dag$ \,and\, Abdelaziz Rhandi$^\dag$}
\date{}

\begin{document}

%Title and Author Information
\maketitle
\vspace{-2em}

\vspace{1em}
\thispagestyle{firstpage}

\begin{abstract}
We investigate the null and approximate controllability of coupled
linear forward stochastic reaction--convection--diffusion systems under
suitable cascade coupling conditions. The model consists of two forward
stochastic parabolic equations governed by general second-order
differential operators with time-, space-, and random-dependent
coefficients. We consider a localized control acting on the drift term
of the first equation, together with controls acting on the diffusion
terms. By means of a duality argument, the controllability problem is
reduced to an observability problem for the associated adjoint
backward stochastic parabolic system. We establish a new global
Carleman estimate for coupled backward stochastic parabolic systems
whose drift terms belong to a negative Sobolev space. As an application
of this estimate, we derive the required observability properties and
establish the null and approximate controllability of the original
forward system.
\end{abstract}

\maketitle
\smallskip


\section{Introduction and Main Results}

Mathematical control theory is a branch of applied mathematics concerned with
the behavior of dynamical systems and the ways in which external inputs can
influence or steer them. A central objective is to determine how to guide a
system's trajectory toward a prescribed outcome, a property known as
controllability. Comprehensive discussions of controllability for deterministic
control systems can be found in
\cite{coron07, fernandez2006global, lions1972some, Zab95, zuazua2007controllability}.
For some controllability results for stochastic systems, we refer the reader to
\cite{gureSiam07, liu2014global, lu2011some, luZhang22mcrf, lu2021mathematical,
tang2009null} and the references therein.

In this work, we study the null and approximate controllability properties of a
cascade system of two coupled forward stochastic parabolic equations with
reaction, convection, and diffusion operators. In the existing literature,
only a limited number of controllability results are available for coupled
stochastic parabolic equations. We refer to
\cite{Preprielgr24, Fadilicasc, LiQi} for coupled backward stochastic
parabolic equations, and to
\cite{liu2014global, obeman255, yansun2011} for coupled forward--backward
stochastic parabolic equations. Results concerning coupled fourth-order
backward stochastic parabolic equations can be found in
\cite{wangNull24}. Furthermore, \cite{liu14couplfor} investigates coupled
forward fractional stochastic parabolic equations, while
\cite{LiuuLiuX} considers a class of forward stochastic parabolic systems.
In these latter two works
\cite{liu14couplfor, LiuuLiuX}, the coefficients are independent of the
spatial variable, and controllability is achieved using the
Lebeau--Robbiano strategy with a single control force. In contrast, in the
present work, we consider the more general setting in which the coefficients
depend on time, space, and random variables.

The primary tool in our analysis is a global Carleman estimate for coupled
backward stochastic parabolic systems with zero- and first-order terms whose
coefficients belong to $L^\infty$, subject to Dirichlet boundary conditions. Introduced by Carleman in 1939 \cite{Carl39} to establish uniqueness for second-order elliptic partial differential equations, Carleman estimates have become fundamental tools in the analysis of deterministic and stochastic partial differential equations. They play a central role in the study of a wide range of problems, including controllability, inverse problems, multi-objective control, and insensitizing control. For comprehensive accounts
and further applications of Carleman estimates in deterministic and stochastic
settings, we refer the reader to
\cite{fernandez2006global, lu2021mathematical, tang2009null, Yamam2009invePrb}.

Systems with coupled interactions arise in numerous applications, including population dynamics, chemical reaction processes, and materials with memory effects. The coupling terms describe the interactions between the different components of the system; see \cite{allafrasl22} and the references therein. Meanwhile, stochastic perturbations account for the intrinsic uncertainties present in realistic environments.

Although the controllability of deterministic coupled reaction--convection--
diffusion systems has been extensively studied (see, e.g.,
\cite{ferluzua16} and the references therein), the corresponding stochastic
systems remain considerably less understood. From a mathematical perspective,
the combination of coupling, convection effects, and randomness introduces
substantial analytical challenges, particularly in the derivation of
observability inequalities for the associated backward stochastic systems.
The presence of martingale terms in the adjoint equations makes the extension
of deterministic controllability techniques to the stochastic framework
highly nontrivial.

The controllability of coupled stochastic parabolic systems has received
increasing attention in recent years. Several results have been obtained for
forward and backward stochastic equations, including systems with fractional
operators, higher-order operators, and coupled forward--backward structures;
see, for instance,
\cite{liu2014global, liu14couplfor, LiuuLiuX, wangNull24} and the references
therein. However, many important questions remain open for coupled forward
stochastic reaction--convection--diffusion systems with spatially dependent
coefficients and general coupling structures.

The aim of the present work is to contribute to this research direction by establishing null and approximate controllability results for a  coupled forward stochastic reaction--convection--diffusion system. In addition, we provide an explicit estimate of the null control cost in terms of the final time $T$ and the coefficients of the system. Our approach is based on the derivation of an appropriate global Carleman estimate for the associated adjoint stochastic system. This estimate enables us to establish the required observability inequality and, through a duality argument, to deduce the controllability results for the original forward stochastic system. To the best of the authors' knowledge, this is the first work addressing the controllability of the coupled case of stochastic parabolic systems extending the single-equation framework previously studied in the literature; see, for instance, \cite{elgconvec23,liu2014global,tang2009null}.

Let \(G \subset \mathbb{R}^N\) (\(N \ge 1\)) be a nonempty, bounded, open domain with closure \(\overline{G}\) and boundary \(\Gamma := \partial G\) of class \(C^4\), and let \(T>0\) denote a fixed final time. We denote by \(G_0\) a nonempty open subset of \(G\), which serves as the control region, and by \(\chi_{G_0}\) its characteristic function. Throughout the paper, we use the following notations:
\[
Q := (0,T) \times G, 
\qquad 
\overline{Q} := [0,T] \times \overline{G}, 
\qquad 
\Sigma := (0,T) \times \Gamma, 
\qquad 
Q_0 := (0,T) \times G_0.
\]

Let \( (\Omega, \mathcal{F}, \{\mathcal{F}_t\}_{t \ge 0}, \mathbb{P}) \) be a fixed complete filtered probability space on which a one-dimensional standard Brownian motion 
\( W(\cdot)\) is defined. 
Here, \( \{\mathcal{F}_t\}_{t \ge 0} \) denotes the natural filtration generated by \( W(\cdot) \), augmented by all \( \mathbb{P} \)-null sets in \( \mathcal{F} \). Let \( \mathcal{X} \) be a Banach space. We denote by \( C([0,T]; \mathcal{X}) \) the Banach space of all continuous \( \mathcal{X} \)-valued functions $\varphi$ on \( [0,T] \), endowed with the norm
$\|\varphi\|_{C([0,T];\mathcal{X})}
:=
\max_{0 \le t \le T} \|\varphi(t)\|_{\mathcal{X}}$. For a sub-\(\sigma\)-algebra \( \mathcal{G} \subset \mathcal{F} \), we define 
\( L^2_{\mathcal{G}}(\Omega; \mathcal{X}) \) as the space of all 
\( \mathcal{X} \)-valued, \( \mathcal{G} \)-measurable random variables \( X \) such that
\(
\mathbb{E}\big[ \|X\|_{\mathcal{X}}^2 \big] < \infty,
\) equipped with the norm
\[
\|X\|_{L^2_{\mathcal{G}}(\Omega;\mathcal{X})}
:=
\big( \mathbb{E}\|X\|_{\mathcal{X}}^2 \big)^{1/2}.
\]
We denote by \( L^2_{\mathcal{F}}(0,T; \mathcal{X}) \) the Banach space of all 
\( \mathcal{X} \)-valued, \( \{\mathcal{F}_t\}_{t \ge 0} \)-adapted stochastic processes 
\( X(\cdot) \) such that
\(
\mathbb{E} \int_0^T \|X(t)\|_{\mathcal{X}}^2 \, dt < \infty,
\) endowed with the norm
\[
\|X\|_{L^2_{\mathcal{F}}(0,T;\mathcal{X})}
:=
\left( \mathbb{E} \int_0^T \|X(t)\|_{\mathcal{X}}^2 \, dt \right)^{1/2}.
\]
Similarly, \( L^\infty_{\mathcal{F}}(0,T; \mathcal{X}) \) denotes the Banach space of all 
\( \mathcal{X} \)-valued, \( \{\mathcal{F}_t\}_{t \ge 0} \)-adapted processes 
\( X(\cdot) \) that are essentially bounded, i.e.,
\(
\operatorname*{ess\,sup}_{(t,\omega)\in(0,T)\times\Omega}
\|X(t,\omega)\|_{\mathcal{X}} < \infty,
\) equipped with the norm
\[
\|X\|_{L^\infty_{\mathcal{F}}(0,T;\mathcal{X})}
:=
\operatorname*{ess\,sup}_{(t,\omega)\in(0,T)\times\Omega}
\|X(t,\omega)\|_{\mathcal{X}}.
\]
Moreover, we define 
\( L^2_{\mathcal{F}}(\Omega; C([0,T]; \mathcal{X})) \) 
as the Banach space of all 
\( \mathcal{X} \)-valued, \( \{\mathcal{F}_t\}_{t \ge 0} \)-adapted continuous processes 
\( X(\cdot) \) satisfying
\(
\mathbb{E} \left[ 
\max_{0 \le t \le T} \|X(t)\|_{\mathcal{X}}^2 
\right] < \infty,
\) with norm
\[
\|X\|_{L^2_{\mathcal{F}}(\Omega;C([0,T];\mathcal{X}))}
:=
\left(
\mathbb{E} \Big[
\max_{0 \le t \le T} \|X(t)\|_{\mathcal{X}}^2
\Big]
\right)^{1/2}.
\]
Similarly, for any positive integer \( n \), one can define 
\( L^\infty_{\mathcal{F}}(\Omega; C^n([0,T]; \mathcal{X})) \). For further details on these stochastic notions, we refer the reader to 
\cite{luZhang22mcrf} and \cite[Chapter~2]{lu2021mathematical}.

In this paper, we use the notation 
\(
\frac{\partial \mathbf{y}}{\partial x_i}
\)
to denote the first-order partial derivative of a function \( \mathbf{y} \) with respect to the variable \( x_i \), where \( x_i \) is the \( i \)-th coordinate of a generic point 
\( x = (x_1, \dots, x_N) \in \mathbb{R}^N \). 
Similarly, 
\(
\frac{\partial^2 \mathbf{y}}{\partial x_i \partial x_j}
\)
denotes the second-order partial derivative with respect to \( x_i \) and \( x_j \). For any \( x, y \in \mathbb{R}^N \), we denote by \( x \cdot y \) their Euclidean inner product, and by \( |x| \) the associated Euclidean norm.

Throughout this work, we assume that the coefficients 
\(
\tau_{ij}^m : Q \times \Omega \to \mathbb{R}
\)
(\( i,j = 1,2,\dots,N \), \( m = 1,2 \)) satisfy the following conditions:
\begin{enumerate}[(1)]
    \item 
    \(
    \tau_{ij}^m \in 
    L^\infty_{\mathcal F}\!\big(\Omega; C^1([0,T]; W^{2,\infty}(G))\big)
    \)
    and
    \(
    \tau_{ij}^m = \tau_{ji}^m
    \)
    for all \( 1 \le i,j \le N \) and \( m=1,2 \).

    \item 
    There exists a constant \( \tau_0 > 0 \) such that
    \begin{equation}\label{assmponalpha}
        \sum_{i,j=1}^N 
        \tau_{ij}^m(t,x,\omega)\,\kappa_i \kappa_j
        \ge 
        \tau_0 |\kappa|^2,
        \qquad 
        \text{for all } 
        (t,x,\omega,\kappa) \in 
        Q \times \Omega \times \mathbb{R}^N,
        \quad m=1,2,
    \end{equation}
    where \( x = (x_1,\dots,x_N) \) and 
    \( \kappa = (\kappa_1,\dots,\kappa_N) \).
\end{enumerate}
For \( m = 1,2 \), we define the second-order differential operators 
\( L_m(t) \) by
\begin{equation}\label{definofopeLi}
L_m(t)\mathbf{y}
:=
\sum_{i,j=1}^N
\frac{\partial}{\partial x_i}
\left(
\tau_{ij}^m(t,x,\omega)
\frac{\partial \mathbf{y}}{\partial x_j}
\right).
\end{equation}

The main objective of this paper is to investigate the null and approximate 
controllability of the following coupled forward stochastic 
reaction--convection--diffusion parabolic system:
\begin{equation}\label{ass15cont}
\begin{cases}
\begin{array}{ll}
dy - L_1(t)y \,dt 
= \Big[\xi_1 + a_{11} y + a_{12} z 
+ B_{11}\cdot\nabla y + B_{12}\cdot\nabla z 
+ f \chi_{G_0} \Big] dt 
+ g_1 \, dW(t) 
& \textnormal{in } Q, \\[0.3em]

dz - L_2(t)z \,dt 
= \Big[\xi_2 + a_{21} y + a_{22} z  
+ B_{21}\cdot\nabla y + B_{22}\cdot\nabla z\Big] dt 
+ g_2 \, dW(t) 
& \textnormal{in } Q, \\[0.3em]

y = z = 0 
& \textnormal{on } \Sigma, \\[0.3em]

y(0) = y_0, \qquad z(0) = z_0 
& \textnormal{in } G,
\end{array}
\end{cases}
\end{equation}
where $(y,z)$ is the state variable, $(y_0,z_0) \in L^2_{\mathcal{F}_0}(\Omega; L^2(G;\mathbb{R}^2))$ is the initial state, and $\xi_1, \xi_2 \in L^2_{\mathcal F}(0,T; L^2(G))$ are some given source terms.
The coefficients fulfill
\[
a_{ij} \in L^\infty_{\mathcal F}(0,T; L^\infty(G)), 
\qquad
B_{ij} \in L^\infty_{\mathcal F}(0,T; L^\infty(G;\mathbb{R}^N)),
\]
for \( i,j \in\{ 1,2\} \).
The control functions satisfy
\[
(f,g_1,g_2) 
\in 
L^2_{\mathcal F}(0,T; L^2(G_0)) 
\times 
\big(L^2_{\mathcal F}(0,T; L^2(G))\big)^2.
\]
In this work, under suitable assumptions on the coupling terms 
$a_{21}$ and $B_{21}$, we study the following controllability problems:
\begin{itemize}
\item First, we prove the \emph{null controllability} of \eqref{ass15cont}, namely, that there exists a constant $C_T>0$ such that, for any initial data $(y_0, z_0)\in L^2_{\mathcal{F}_0}(\Omega; L^2(G;\mathbb{R}^2))$ and source terms $(\xi_1,\xi_2)\in\Xi$ (a subspace of $L^2_{\mathcal F}(0,T; L^2(G;\mathbb{R}^2))$), one can find controls $(f,g_1,g_2)\in L^2_{\mathcal F}(0,T; L^2(G_0)) \times \big(L^2_{\mathcal F}(0,T; L^2(G))\big)^2$ such that the corresponding solution $(y, z)$ satisfies
\[
y(T,\cdot)=z(T,\cdot)=0 \quad \text{in } G,\ \mathbb{P}\text{-a.s.},
\]
and
\[
\|(f,g_1,g_2)\|^2_{L^2_{\mathcal F}(0,T; L^2(G_0)) \times \big(L^2_{\mathcal F}(0,T; L^2(G))\big)^2}
\le C_T\Big(\|(y_0, z_0)\|^2_{L^2_{\mathcal{F}_0}(\Omega; L^2(G;\mathbb{R}^2))}+\|(\xi_1,\xi_2)\|^2_\Xi\Big).
\]
\item Next, we establish the \emph{approximate controllability} of \eqref{ass15cont}, namely, that for any initial data $(y_0, z_0)\in L^2_{\mathcal{F}_0}(\Omega; L^2(G;\mathbb{R}^2))$, target states $(y_T, z_T)\in L^2_{\mathcal{F}_T}(\Omega; L^2(G;\mathbb{R}^2))$, and source terms $(\xi_1,\xi_2)\in L^2_{\mathcal F}(0,T; L^2(G;\mathbb{R}^2))$, one can find controls $(f,g_1,g_2)\in L^2_{\mathcal F}(0,T; L^2(G_0)) \times \big(L^2_{\mathcal F}(0,T; L^2(G))\big)^2$ such that the corresponding solution $(y, z)$ satisfies
\[
(y(T,\cdot), z(T,\cdot)) \ \text{is arbitrarily close to} \ (y_T, z_T) \ \text{in }\, L^2_{\mathcal{F}_T}(\Omega; L^2(G;\mathbb{R}^2)).
\]
\end{itemize}
 
By the classical duality argument (see, e.g.,
\cite[Chapter 7]{lu2021mathematical}), it is well known that the null and
approximate controllability of \eqref{ass15cont} can be reduced, respectively,
to an observability inequality and a uniqueness property for the following
coupled adjoint backward stochastic system:
\begin{equation}\label{adback4.77firstcontro}
\begin{cases}
\begin{array}{ll}
dh + L_1(t) h \, dt 
= \Big[-a_{11} h - a_{21} k  + \operatorname{div}(h B_{11} + k B_{21}) \Big] dt 
+ H \, dW(t) 
& \text{in } Q, \\[0.3em]

dk + L_2(t) k \, dt 
= \Big[-a_{12} h - a_{22} k  + \operatorname{div}(h B_{12} + k B_{22}) \Big] dt 
+ K \, dW(t) 
& \text{in } Q, \\[0.3em]

h = k = 0 
& \text{on } \Sigma, \\[0.3em]

h(T) = h_T, \quad k(T) = k_T 
& \text{in } G,
\end{array}
\end{cases}
\end{equation}
where \( (h,k;H,K) \) denotes the state variable,
\((h_T,k_T)\in L^2_{\mathcal{F}_T}(\Omega;L^2(G;\mathbb{R}^2))\) is the
 terminal state.

Concerning the well-posedness of systems \eqref{ass15cont} and
\eqref{adback4.77firstcontro}, we consider the Gelfand triple
\[
H_0^1(G;\mathbb{R}^2)
\subset
L^2(G;\mathbb{R}^2)
\subset
H^{-1}(G;\mathbb{R}^2).
\]

For the forward system \eqref{ass15cont}, we introduce the state vector
\[
U(t):=(y(t),z(t))\in L^2(G;\mathbb{R}^2),
\]
and rewrite the system in the abstract form
\begin{equation}\label{abstract-SPDE}
\begin{cases}
dU(t)=\big(\mathcal A(t)U(t)+F(t)\big)\,dt
+\widetilde F(t)\,dW(t),
& t\in(0,T],\\[0.3em]
U(0)=U_0,
\end{cases}
\end{equation}
where
\[
F(t)=
\begin{pmatrix}
\xi_1+f\chi_{G_0}\\
\xi_2
\end{pmatrix},
\qquad
\widetilde F(t)=
\begin{pmatrix}
g_1\\
g_2
\end{pmatrix},
\]
and the operator
\[
\mathcal A(t):H_0^1(G;\mathbb{R}^2)
\longrightarrow H^{-1}(G;\mathbb{R}^2)
\]
is defined by
\[
\mathcal A(t)
\begin{pmatrix}
y\\ z
\end{pmatrix}
=
\begin{pmatrix}
L_1(t)y+a_{11}y+a_{12}z
+B_{11}\cdot\nabla y+B_{12}\cdot\nabla z
\\[1mm]
L_2(t)z+a_{21}y+a_{22}z
+B_{21}\cdot\nabla y+B_{22}\cdot\nabla z
\end{pmatrix}.
\]
It is straightforward to verify that $\mathcal{A}(t)$ is linear and bounded from
$H_0^1(G;\mathbb{R}^2)$ into $H^{-1}(G;\mathbb{R}^2)$, and that it satisfies the stochastic parabolicity condition in \cite[Theorem~4.4.3]{lotorovozski}. Moreover, the forcing terms $F(t)$ and $\widetilde{F}(t)$ are $\{\mathcal{F}_t\}$-adapted processes with values in $L^2(G;\mathbb{R}^2)$ and $\mathcal{L}(\mathbb{R}^2;L^2(G;\mathbb{R}^2))$, respectively. Consequently, all the assumptions of \cite[Theorem~4.4.3]{lotorovozski} are satisfied, and then system \eqref{abstract-SPDE} (and hence \eqref{ass15cont}) is well posed. More precisely, for any initial datum $U_0\in L^2_{\mathcal{F}_0}(\Omega;L^2(G;\mathbb{R}^2))$
and forcing terms $F\in L^2_{\mathcal{F}}(0,T;L^2(G;\mathbb{R}^2)),$ and $\widetilde{F}\in L^2_{\mathcal{F}}(0,T;\mathcal{L}(\mathbb{R}^2;L^2(G;\mathbb{R}^2))),$ there exists a unique weak solution
\[
U=(y,z)\in
L^2_{\mathcal F}\big(\Omega;C([0,T];L^2(G;\mathbb R^2))\big)
\cap
L^2_{\mathcal F}(0,T;H_0^1(G;\mathbb R^2)).
\]
Furthermore, the following estimate holds:
\[
\begin{aligned}
&\|U\|_{L^2_{\mathcal F}(\Omega;C([0,T];L^2(G;\mathbb R^2)))}
+\|U\|_{L^2_{\mathcal F}(0,T;H_0^1(G;\mathbb R^2))}
\\
&\qquad\leq C_T\Big(
\|U_0\|_{L^2_{\mathcal F_0}(\Omega;L^2(G;\mathbb R^2))}
+\|F\|_{L^2_{\mathcal F}(0,T;L^2(G;\mathbb R^2))}
+\|\widetilde F\|_{L^2_{\mathcal F}
(0,T;L^2(G;\mathbb R^2)))}
\Big),
\end{aligned}
\]
where $C_T>0$ is independent of the data.

On the other hand, the adjoint system
\eqref{adback4.77firstcontro} can be written in the abstract backward
stochastic evolution form
\begin{equation}\label{esuq11.lkback}
\begin{cases}
d\mathbf h(t)
=
-\mathcal A^*(t)\mathbf h(t)\,dt
+\mathbf R(t)\,dW(t),
& t\in[0,T),\\[0.3em]
\mathbf h(T)=\mathbf h_T
\in L^2_{\mathcal F_T}(\Omega;L^2(G;\mathbb R^2)),
\end{cases}
\end{equation}
where $\mathcal A^*(t):H_0^1(G;\mathbb R^2)
\rightarrow H^{-1}(G;\mathbb R^2)$ denotes the adjoint operator of
$\mathcal A(t)$.

The adjoint variables and the martingale term are defined by
\[
\mathbf h(t):=
\begin{pmatrix}
h(t)\\[1mm]
k(t)
\end{pmatrix},
\qquad
\mathbf R(t):=
\begin{pmatrix}
H(t)\\[1mm]
K(t)
\end{pmatrix}.
\]

Under the same assumptions ensuring the well-posedness of the forward
system, the operator $\mathcal A^*(t)$ satisfies the standard assumptions
for backward stochastic parabolic equations. Hence, by the classical theory
of backward stochastic evolution equations (see, e.g.,
\cite{hupeng91}), system \eqref{esuq11.lkback} (and hence \eqref{adback4.77firstcontro}) admits a unique weak solution
\[
(\mathbf h,\mathbf R)
\in
\Big(
L^2_{\mathcal F}\big(\Omega;C([0,T];L^2(G;\mathbb R^2))\big)
\cap
L^2_{\mathcal F}(0,T;H_0^1(G;\mathbb R^2))
\Big)
\times
L^2_{\mathcal F}
\big(0,T;L^2(G;\mathbb R^2))\big).
\]
Moreover, there exists a constant $C_T>0$, independent of the terminal datum
$\mathbf h_T$, such that
\[
\begin{aligned}
&\|\mathbf h\|_{L^2_{\mathcal F}
(\Omega;C([0,T];L^2(G;\mathbb R^2)))}
+
\|\mathbf h\|_{L^2_{\mathcal F}
(0,T;H_0^1(G;\mathbb R^2))}
\\
&\qquad
+
\|\mathbf R\|_{L^2_{\mathcal F}
(0,T;L^2(G;\mathbb R^2)))}
\\
&\leq
C_T\,
\|\mathbf h_T\|_{L^2_{\mathcal F_T}
(\Omega;L^2(G;\mathbb R^2))}.
\end{aligned}
\]

In order to establish the null and approximate controllability results for system
\eqref{ass15cont}, we assume that the
coupling coefficients $a_{21}$ and $B_{21}$ satisfy the following \emph{cascade structure}
conditions:

\begin{itemize}

\item[(i)] \textbf{Zero-order coupling coefficient $a_{21}$.}
\begin{equation}\label{assump1.3}
a_{21}\geq a_0>0
\quad \text{or}\quad
-a_{21}\geq a_0>0,
\quad
\text{in }(0,T)\times\widetilde G_0,
\quad \mathbb P\text{-a.s.},
\end{equation}
for a nonempty open subset $\widetilde G_0\Subset G_0$.

\item[(ii)] \textbf{First-order coupling coefficient $B_{21}$.}
\begin{equation}\label{assump1.3B21}
B_{21}=0,
\quad
\text{in }Q,
\quad \mathbb P\text{-a.s.}
\end{equation}

\end{itemize}

The above cascade structure conditions (i)--(ii) ensure that the localized
control $f$ acts directly on the first component $y$ of system
\eqref{ass15cont} and that its influence is transferred to the second
component $z$ through the coupling terms. Condition (i) guarantees a
non-degenerate zero-order coupling $a_{21}$ from $y$ to $z$ in a subdomain,
which plays a fundamental role in recovering observability of the second
component from the first one. Condition (ii) eliminates the first-order
coupling term $B_{21}$, which allows the Carleman estimate to be closed
without additional regularity assumptions on the coupling coefficients. Together, these assumptions create a cascade structure: the control acting on
the first equation propagates through the coupling mechanism and provides
information on the second equation. This structure is essential for deriving
the observability inequality of the adjoint system and, consequently, for
establishing the controllability results. Further discussions on these
assumptions are given in Remark~\ref{rmkk13}.

Throughout this paper, $C$ denotes a generic positive constant depending only on $G$, $G_0$, $\widetilde{G}_0$, and the coefficients of the system \eqref{ass15cont}. Its value may change from line to line.

The first main result of this paper is the following null controllability
result for system \eqref{ass15cont}.

\begin{thm}\label{thmm1.3insfirst}
Assume that \eqref{assump1.3} and \eqref{assump1.3B21} hold. Then, for any initial states $y_0,z_0\in L^2_{\mathcal F_0}(\Omega;L^2(G))$,
and any source terms $\xi_1,\xi_2\in L^2_{\mathcal F}(0,T;L^2(G))$ satisfying
\begin{align}\label{assonxi1xi2}
\max\bigg\{
\left\|\theta^{-1}\gamma^{-3}\xi_1\right\|_
{L^2_{\mathcal F}(0,T;L^2(G))},
\left\|\theta^{-1}\gamma^{-3/2}\xi_2\right\|_
{L^2_{\mathcal F}(0,T;L^2(G))}
\bigg\}<\infty,
\end{align}
there exists a control triple
\[
(f,g_1,g_2)
\in
L^2_{\mathcal F}(0,T;L^2(G_0))
\times
\big(L^2_{\mathcal F}(0,T;L^2(G))\big)^2
\]
such that the corresponding solution $(y,z)$ of \eqref{ass15cont}
satisfies
\[
y(T,\cdot)=z(T,\cdot)=0
\quad\textnormal{in }G,\qquad
\mathbb P\textnormal{-a.s.}
\]
Moreover, there exists a constant $C>0$ such that
\begin{align*}
&\|f\|_{L^2_{\mathcal F}(0,T;L^2(G_0))}
+\|g_1\|_{L^2_{\mathcal F}(0,T;L^2(G))}
+\|g_2\|_{L^2_{\mathcal F}(0,T;L^2(G))}
\\
&\leq
\sqrt{\exp(CM)}
\bigg(
\|y_0\|_{L^2_{\mathcal F_0}(\Omega;L^2(G))}
+\|z_0\|_{L^2_{\mathcal F_0}(\Omega;L^2(G))}
\\
&\qquad\qquad
+\|\theta^{-1}\gamma^{-3}\xi_1\|_
{L^2_{\mathcal F}(0,T;L^2(G))}
+\|\theta^{-1}\gamma^{-3/2}\xi_2\|_
{L^2_{\mathcal F}(0,T;L^2(G))}
\bigg),
\end{align*}
where $\theta$ and $\gamma$ are the weight functions defined in
\eqref{2.2012}, and
\begin{align*}
    \begin{aligned}
M&=1+T+\frac{1}{T}+\|a_{11}\|_\infty^{\frac{2}{3}}+\|B_{11}\|^2_\infty +\|a_{22}\|_\infty^{\frac{2}{3}}+\|B_{22}\|^2_\infty +\|a_{12}\|_\infty^{\frac{1}{3}}+\|B_{12}\|_\infty^{\frac{1}{2}}\\
&\quad+T(\|a_{11}\|_\infty+\|a_{12}\|_\infty+\|a_{22}\|_\infty+\|B_{11}\|^2_\infty+\|B_{12}\|^2_\infty+\|B_{22}\|^2_\infty).
    \end{aligned}
\end{align*}
\end{thm}

The second main result is the following approximate controllability for
system \eqref{ass15cont}.

\begin{thm}\label{thm1.2approx}
Assume that 
\eqref{assump1.3} and \eqref{assump1.3B21} hold. Then, for any $\xi_1,\xi_2\in L^2_{\mathcal F}(0,T;L^2(G))$, any initial states $y_0,z_0\in L^2_{\mathcal F_0}(\Omega;L^2(G))$, any terminal states
$y_T,z_T\in L^2_{\mathcal F_T}(\Omega;L^2(G))$, and for every $\varepsilon>0$, there exist controls
\[
(f,g_1,g_2)\in
L^2_{\mathcal F}(0,T;L^2(G_0))
\times
\big(L^2_{\mathcal F}(0,T;L^2(G))\big)^2,
\]
such that the corresponding solution $(y,z)$ of
\eqref{ass15cont} satisfies
\[
\mathbb E\|y(T)-y_T\|_{L^2(G)}^2\leq\varepsilon\quad\textnormal{and}\quad
\mathbb E\|z(T)-z_T\|_{L^2(G)}^2\leq\varepsilon .
\]
\end{thm}

Unlike the controllability of coupled deterministic parabolic systems
(see, e.g., \cite{surveyAmmarKBGT, GonBurTer}), the system
\eqref{ass15cont} considered in this work involves three control inputs:
\( f \), which acts locally on the drift term of the first equation, and two
additional controls, \( g_1 \) and \( g_2 \), which influence the diffusion
terms globally. The introduction of the controls \( g_1 \) and \( g_2 \) is
motivated by the stochastic nature of the equations, where the Brownian
perturbations generate additional martingale terms in the associated adjoint
system. Moreover, these controls are closely related to the structure of the
backward stochastic equations arising in the dual formulation of the
controllability problem. Indeed, the adjoint problem is a coupled backward stochastic parabolic system
involving both the state variables and the martingale correction terms
\(H\) and \(K\). These additional stochastic components ensure the
well-posedness of the adjoint system but also create substantial difficulties
in the derivation of suitable Carleman estimates. Such difficulties do not
appear in the same form for deterministic parabolic systems or scalar
stochastic equations; see, e.g.,
\cite{elgconvec23,tang2009null}.

In this paper, we impose the cascade assumptions \eqref{assump1.3} and \eqref{assump1.3B21} on the coupling terms \(a_{21}\) and \(B_{21}\), respectively.
 These assumptions guarantee that the
localized control acting on the first component is transferred to the second
component through the coupling structure. The positivity condition on
\(a_{21}\) in a subdomain plays a fundamental role in the proof of the
Carleman estimate and the resulting observability inequality. Moreover, the
condition imposed on \(B_{21}\) allows us to handle the first-order coupling
terms arising in the weighted energy estimates. Without these structural
conditions, the coupling mechanism may not be sufficient to obtain the
desired controllability properties. To illustrate the importance of these
assumptions, we provide the following example:
\begin{itemize}
\item Let \( \xi_i \equiv B_{ij}\equiv 0 \), $a_{ij}\equiv0$ for $(i,j)\neq(1,2)$, \( \tau_{ij}^m = \delta_{ij} \) (where \( \delta_{ij} \) denotes the Kronecker delta with \( 1 \leq i, j \leq N \)), then the system \eqref{ass15cont} reduces to the following coupled stochastic system:
\begin{equation}\label{ass15contsec}
    \begin{cases}
    \begin{array}{ll}
    dy - \Delta y \, dt = (a_{12}z+f \chi_{G_0} ) \, dt + g_1 \, dW(t) & \text{in } Q, \\
    dz - \Delta z \, dt = g_2 \, dW(t) & \text{in } Q, \\
    y = z = 0 & \text{on } \Sigma, \\
    y(0) = y_0, \quad z(0) = z_0 & \text{in } G.
    \end{array}
    \end{cases}
\end{equation}

We claim that the system \eqref{ass15contsec} is neither null controllable nor approximately controllable at time $T$. To demonstrate this, we construct a particular solution to the adjoint equation that does not satisfy the appropriate observability property. Indeed, the adjoint equation associated with \eqref{ass15contsec} is the following backward system:
\begin{equation}\label{adbackadjointex}
\begin{cases}
\begin{array}{ll}
dh + \Delta h \, dt = H \, dW(t) & \text{in } Q, \\
dk + \Delta k \, dt = -a_{12}h\,dt + K \, dW(t) & \text{in } Q, \\
h = k = 0 & \text{on } \Sigma, \\
h(T) = h_T, \quad k(T) = k_T & \text{in } G.
\end{array}
\end{cases}
\end{equation}
By a classical duality argument, system \eqref{ass15contsec} is null controllable at time \( T \) if and only if the solutions of \eqref{adbackadjointex} satisfy the observability inequality
\begin{align}\label{obseineex}
\mathbb{E}\int_G \left(|h(0)|^2 + |k(0)|^2\right)\,dx
\leq C_T\, \mathbb{E} \iint_Q \left(  h^2\chi_{G_0} + H^2 + K^2 \right) \, dx \, dt,
\end{align}
for some positive constant \( C_T \). Moreover, system \eqref{ass15contsec} is approximately controllable at time \( T \) if and only if the solutions of \eqref{adbackadjointex} satisfy the unique continuation property
\begin{align}\label{uniqprot}
 (h, H, K) = (0, 0, 0) \quad \text{in } Q_0 \times Q^2, \quad \mathbb{P}\text{-a.s.}
 \quad \Longrightarrow \quad
 h_T = k_T = 0 \quad \text{in } G, \quad \mathbb{P}\text{-a.s.}
\end{align}
Let \( s_1 \) be the first eigenvalue of \( -\Delta \) (the Dirichlet Laplacian) and \( e_1\in H^2(G)\cap H^1_0(G) \) the associated eigenfunction such that \( \|e_1\|_{L^2(G)} = 1 \). Define the function \( p(t) := \nu(t)e_1 \), where \( \nu \in C^1([0,T]) \) is the solution of the following differential equation
\begin{align*}
    \begin{cases}
   \nu' - s_1 \nu = 0, \quad \text{in } (0,T), \\
   \nu(0) = 1.     
    \end{cases}
\end{align*}
Notice that \( p \) solves the following backward heat equation
\begin{align}\label{eqbacksybyg}
\begin{cases}
\begin{array}{ll}
\partial_t p + \Delta p =  0 & \text{in } Q, \\p = 0 & \text{on } \Sigma.
\end{array}
\end{cases}
\end{align}
By the well-posedness result for backward stochastic equations, it follows that the solution to \eqref{adbackadjointex} with \( h_T = 0 \) and \( k_T = \nu(T)e_1 \) is given by \( (h, k; H, K) \equiv (0, p; 0, 0) \). It is then straightforward to verify that this solution does not satisfy the observability inequality \eqref{obseineex} nor the unique continuation property \eqref{uniqprot}. Consequently, the system \eqref{ass15contsec} is neither null controllable nor approximately controllable. This  example illustrates that the additional control \( g_2 \) is insufficient to drive the second component \( z \) of the system \eqref{ass15contsec} to zero at the final time \( T \). In the absence of the condition \eqref{assump1.3}, the influence of \( g_2 \) on the diffusion term in the second equation of \eqref{ass15contsec} is significantly restricted. Consequently, this limitation prevents the establishment of the controllability results and highlights the necessity of the assumption \eqref{assump1.3} for the control of coupled forward stochastic parabolic systems.
\end{itemize}
Notice that no additional control acting on the drift term of the second
equation in \eqref{ass15cont} is required, thanks to assumption
\eqref{assump1.3} on the coupling term $a_{21}y$, which provides an indirect
control mechanism for the drift part of the second equation. On the other
hand, using only the control $f$, without the additional controls $g_1$ and
$g_2$, to address controllability problems for general forward stochastic
parabolic systems whose coefficients depend on time, space, and random
variables remains an open problem. For further details on partial
controllability results with a single control for a class of stochastic heat
equations, we refer to \cite{withouextra, lu2011some, observineqback}. For
results concerning coupled forward stochastic parabolic equations, see
\cite{liu14couplfor, LiuuLiuX}. In all these works, only space-independent zero-order terms are considered, and the main analytical tool employed is the spectral method based on the Lebeau--Robbiano strategy. In the present work, we investigate coupled forward stochastic parabolic systems involving zero-, first-, and second-order terms with coefficients that also depend on the space variable.

\begin{rmk}\label{rmkk13}
The assumption \eqref{assump1.3B21} plays an essential role in establishing the Carleman estimate \eqref{carlestcasca1.5sec}. Indeed, in Step~2 of the proof, the term
\[
\lambda^3 \, \mathbb{E} \iint_Q \gamma^3 k \, B_{21}\cdot\nabla(\theta^2\rho_1 k)\,dx\,dt
\]
vanishes in the identity \eqref{forB21} due to assumption \eqref{assump1.3B21}. Without this condition, an application of Young's inequality would produce weighted terms involving higher powers of the Carleman parameter \(\lambda\) than those available on the left-hand side of \eqref{carlestcasca1.5sec}. More precisely, the left-hand side only provides control of the weighted quantities
\[
\lambda^3 \, \mathbb{E} \iint_Q \theta^2\gamma^3 k^2\,dx\,dt,
\qquad
\lambda \, \mathbb{E} \iint_Q \theta^2\gamma |\nabla k|^2\,dx\,dt,
\]
which are not sufficient to absorb the additional contributions generated by the first-order coupling term in the right-hand side. Therefore, treating the case \(B_{21}\neq0\) requires a different Carleman argument and is not covered by the present approach.
\end{rmk}
\begin{rmk}
The approach developed in this paper can be extended to investigate the controllability of a system of $n$ forward stochastic parabolic equations by introducing a localized control in the first equation together with $n$ controls $(g_i)_{i=1,\ldots,n}$ acting in the diffusion terms, under cascade coupling conditions similar to \eqref{assump1.3} and \eqref{assump1.3B21}. Furthermore, since the controls $g_1$ and $g_2$ act globally in the diffusion terms of \eqref{ass15cont}, it is reasonable to expect that, once the controllability of the random parabolic system associated with \eqref{ass15cont} is established by means of \emph{adapted controls} $(f,g_1,g_2)$, the controllability of more general stochastic systems involving zero- and first-order terms in the diffusion parts can be established as a direct consequence.
\end{rmk}
\begin{rmk}\label{rmk:2D-BM}
The controllability formulation of \eqref{ass15cont} is given for a
one-dimensional Brownian motion \(W\), which is sufficient for the
controllability results established in this paper. We briefly indicate
the corresponding modifications in the case of a two-dimensional
Brownian motion \(W=(W^1,W^2)\), with the filtration generated by \(W\).
In this setting, the martingale representation theorem requires both
Brownian directions in each component of the backward adjoint system.
Thus, the stochastic terms \(H\,dW\) and \(K\,dW\) in
\eqref{adback4.77firstcontro} should be replaced by
\[
H_1\,dW^1+H_2\,dW^2,
\qquad
K_1\,dW^1+K_2\,dW^2,
\]
respectively. Accordingly, following the proof strategy of
\eqref{carlthm32f1.5}, the corresponding Carleman estimate for the full
backward system contains the stochastic terms
\[
|H_1|^2+|H_2|^2+|K_1|^2+|K_2|^2
\]
on its right-hand side. The associated It\^o product formulas are
modified according to
\[
d[W^1,W^2]_t=0,
\qquad
d[W^i,W^i]_t=dt,\quad i=1,2.
\]
Consequently, the observability and controllability arguments must be
reformulated consistently to account for these additional martingale
components. In particular, to obtain a controllability result with
the full four-component stochastic observability structure, the
forward system would naturally have to be equipped with four diffusion
controls, one for each Brownian direction in each equation. More
precisely, \eqref{ass15cont} would take the form
\begin{equation*}
\begin{cases}
\begin{array}{ll}
dy - L_1(t)y\,dt
= \Big[\xi_1+a_{11}y+a_{12}z
+B_{11}\cdot\nabla y+B_{12}\cdot\nabla z
+f\chi_{G_0}\Big]dt\\
\hspace{3cm} +g_{11}\,dW^1(t)+g_{12}\,dW^2(t),
& \textnormal{in }Q,\\[0.3em]

dz - L_2(t)z\,dt
= \Big[\xi_2+a_{21}y+a_{22}z
+B_{21}\cdot\nabla y+B_{22}\cdot\nabla z\Big]dt\\
\hspace{3cm}
+g_{21}\,dW^1(t)+g_{22}\,dW^2(t),
& \textnormal{in }Q,\\[0.3em]

y=z=0,
& \textnormal{on }\Sigma,\\[0.3em]

y(0)=y_0,\qquad z(0)=z_0,
& \textnormal{in }G,
\end{array}
\end{cases}
\end{equation*}
where \(g_{11},g_{12},g_{21},g_{22}\) are, respectively, the diffusion
controls paired, through the duality relation, with the stochastic
components \(H_1,H_2,K_1,K_2\) of the backward adjoint system.
\end{rmk}

In what follows, we use the following lemma for integration by parts in space (see \cite{evans}).
\begin{lm}\label{lm1.2}
Let \(F \in L^2(G; \mathbb{R}^N)\). The divergence of \(F\) can be extended as a continuous linear functional on $H^1_0(G)$, defined by
\begin{equation*}
\textnormal{div}(F) : H^1_0(G) \longrightarrow \mathbb{R}, 
\quad u \longmapsto -\int_G F \cdot \nabla u \, dx.
\end{equation*}
Moreover, the following estimate holds:
\begin{align*}
\left|
\langle \textnormal{div}(F), u \rangle_{H^{-1}(G),\, H^1_0(G)}
\right|
&\le C \|F\|_{L^2(G; \mathbb{R}^N)} \|u\|_{H^1_0(G)},
\end{align*}
for all \(u \in H^1_0(G)\), and for some positive constant \(C\).
\end{lm}

We now present the following Itô's formula for stochastic processes in weak form (see, e.g., \cite[Chapter~2]{lu2021mathematical}). This result will be useful for deriving energy estimates and for establishing duality relations between forward and backward systems whose drift terms lie in the negative Sobolev space \(H^{-1}(G)\).
\begin{lm}\label{lm1.1}
Let \(z, y \in L^2_\mathcal{F}(0,T; H^1_0(G))\), \(z_T \in L^2_{\mathcal{F}_T}(\Omega; L^2(G))\), \(y_0 \in L^2_{\mathcal{F}_0}(\Omega; L^2(G))\), \(\phi, \widetilde{\phi} \in L^2_\mathcal{F}(0,T; H^{-1}(G))\), and \(Y, Z \in L^2_\mathcal{F}(0,T; L^2(G))\). Assume that, for all \(t \in [0,T]\), the processes \((z, Z)\) and \(y\) satisfy the following equations in \(G\), respectively:
\begin{align*}
z(t) &= z_T
- \int_t^T \phi(s) \, ds
- \int_t^T Z(s) \, dW(s),
\quad \mathbb{P}\textnormal{-a.s.}, \\
y(t) &= y_0
+ \int_0^t \widetilde{\phi}(s)\, ds
+ \int_0^t Y(s) \, dW(s),
\quad \mathbb{P}\textnormal{-a.s.}
\end{align*}
Then the following assertions hold:
\begin{enumerate}[1.]
\item For any \(t \in [0,T]\), it holds almost surely that
\begin{align*}
\| z(t) \|^2_{L^2(G)}
&= \| z(0) \|^2_{L^2(G)}
+ 2 \int_0^t
\langle \phi(s), z(s)
\rangle_{H^{-1}(G),H^1_0(G)} \, ds \\
&\quad + 2 \int_0^t
\langle z(s), Z(s)
\rangle_{L^2(G)} \, dW(s)
+ \int_0^t
\| Z(s)\|^2_{L^2(G)} \, ds.
\end{align*}

\item For all $t \in [0,T]$, we have almost surely that
\begin{align*}
\langle z(t),y(t)
\rangle_{L^2(G)}
&= \langle z(0),y_0
\rangle_{L^2(G)} + \int_0^t
\langle \phi(s),y(s)
\rangle_{H^{-1}(G),H^1_0(G)} \, ds \\
&\quad + \int_0^t
\langle \widetilde{\phi}(s),z(s)
\rangle_{H^{-1}(G),H^1_0(G)} \, ds + \int_0^t
\langle y(s),Z(s)
\rangle_{L^2(G)} \, dW(s) \\
&\quad + \int_0^t
\langle z(s),Y(s)
\rangle_{L^2(G)} \, dW(s) + \int_0^t
\langle Z(s), Y(s) \rangle_{L^2(G)} \, ds.
\end{align*}
\end{enumerate}
\end{lm}

The rest of this paper is organized as follows. In Section~\ref{sec3}, we derive a suitable Carleman estimate for the adjoint system~\eqref{adback4.77firstcontro}. Section~\ref{sec4} is devoted to the analysis of the associated observability problems, where we establish an observability inequality and a unique continuation property for system~\eqref{adback4.77firstcontro}. Finally, in Section~\ref{sec5}, we prove our main controllability results, namely Theorems~\ref{thmm1.3insfirst} and~\ref{thm1.2approx}.
\section{Global Carleman Estimates}\label{sec3}
In this section, we derive a new global Carleman estimate for the coupled system \eqref{adback4.77firstcontro}. This estimate will be essential for addressing our control problems considered in this paper.

\subsection{Preliminary Tools}
In this subsection, we provide some essential tools for Carleman estimates for backward stochastic parabolic equations. First, we recall the following known result (see \cite{fursikov1996controllability}).
\begin{lm}\label{lmm5.1}
For any nonempty open subset $\mathcal{B}\Subset G$, there exists a function $\psi\in C^4(\overline{G})$ such that
$$
\psi>0\;\,\, \textnormal{in} \,\,G\,;\qquad \psi=0\;\,\,\, \textnormal{on} \,\,\Gamma;\qquad\vert\nabla\psi\vert>0\; \,\,\,\,\textnormal{in}\,\,\overline{G\setminus\mathcal{B}}.
$$
\end{lm}
For any parameters $\lambda, \mu\geq1$, we define the following weight functions
\begin{align}\label{2.2012}
\begin{aligned}
\gamma\equiv\gamma(t)=[t(T-t)]^{-1},\qquad \sigma=\sigma(x)=e^{\mu\psi(x)}-e^{2\mu\|\psi\|_\infty}
\\
\alpha\equiv\alpha(t,x) = \sigma(x)\gamma(t),\qquad\; \theta\equiv\theta(t,x)=e^{\lambda\alpha}.
\end{aligned}
\end{align}
It is easy to check that there exists a constant $C=C(G)>0$ such that for all $(t,x)\in Q$, 
\begin{equation}\label{2.301}
\begin{aligned}
&|\gamma'(t)| \leq C T \gamma^2(t), \qquad &&|\gamma''(t)| \leq C T^2 \gamma^3(t),\\[2mm]
&|\alpha_t(t,x)| \leq C T e^{2\mu\|\psi\|_\infty} \gamma^2(t), \qquad &&|\alpha_{tt}(t,x)| \leq C T^2 e^{2\mu\|\psi\|_\infty} \gamma^3(t),\\[1mm]
&\gamma^{-s}(t) \leq C T^{2s}, \quad\forall s>0.
\end{aligned}
\end{equation}
In what follows, for $d\in\mathbb{R}$, we denote by
\[
\mathcal{I}(d, \cdot) = \lambda^d \mathbb{E} \iint_Q \theta^2 \gamma^d |\cdot|^2 \, dx \, dt + \lambda^{d-2} \mathbb{E} \iint_Q \theta^2 \gamma^{d-2} \vert\nabla \cdot\vert^2 \, dx \, dt.
\] 

Let us consider the following backward stochastic parabolic equation
\begin{equation}\label{eqqgbc}
\begin{cases}
dz+ L_0(t)z \,dt=(F_1+\textnormal{div}(F))\,dt+ Z \,dW(t) & \textnormal{in}\,\,Q,\\ 
z=0 & \textnormal{on}\,\,\Sigma,\\
z(T)=z_T & \textnormal{in}\,\, G,
\end{cases}
\end{equation}
where $z_T \in L^2_{\mathcal{F}_T}(\Omega; L^2(G))$, $F_1 \in L^2_\mathcal{F}(0,T; L^2(G))$, $F \in L^2_\mathcal{F}(0,T; L^2(G;\mathbb{R}^N))$, 
and the operator $L_0(t)$ is defined as
$$L_0(t)z=\displaystyle\sum_{i,j=1}^N \frac{\partial}{\partial x_i}\left(\tau_{ij}^0(t,x,\omega)\frac{\partial z}{\partial x_j}\right),$$
with $\tau_{ij}^0\in L^\infty_\mathcal{F}(\Omega;C^1([0,T];W^{2,\infty}(G)))$, $\tau_{ij}^0=\tau_{ji}^0$ ($1\leq i,j\leq N$),  and there exists a positive constant $\tau^0$ so that 
$$\sum_{i,j=1}^N \tau_{ij}^0(t,x,\omega)\kappa_i\kappa_j\geq \tau^0|\kappa|^2\qquad\textnormal{for any}\quad (t,x,\omega,\kappa)\in Q\times \Omega\times\mathbb{R}^N.$$

Let us recall the following Carleman estimate for equation \eqref{eqqgbc} (see \cite[Theorem 3.1]{elgconvec23}).
\begin{lm}\label{thmm3.1cab}
Let $\mathcal{B}\subset G$ be a nonempty open subset. There exist a large $\mu_0\geq1$ such that for $\mu=\mu_0$, one can find a constant $C>0$ and a large $\lambda_0\geq1$ depending only on $G$, $\mathcal{B}$, $\mu_0$, $\tau_{ij}^0$, and $\tau^0$ such that for all  $F_1\in L^2_\mathcal{F}(0,T;L^2(G))$, $F\in L^2_\mathcal{F}(0,T;L^2(G;\mathbb{R}^N))$ and $z_T\in L^2_{\mathcal{F}_T}(\Omega;L^2(G))$, the solution $(z,Z)$ of equation \eqref{eqqgbc} satisfies 
\begin{align}\label{carback5.8}
\begin{aligned}
\mathcal{I}(3,z)\leq C \bigg[& \lambda^3 \mathbb{E} \int_{0}^T \int_{\mathcal{B}} \theta^2 \gamma^3 z^2 \, dx \, dt+ \mathbb{E}\iint_Q \theta^2 F_1^2 \,dx\,dt\\
&+ \lambda^2\mathbb{E}\iint_Q \theta^2\gamma^2 |F|^2 \,dx\,dt+\lambda^2\mathbb{E}\iint_Q \theta^2\gamma^2 Z^2 \,dx\,dt\bigg], 
\end{aligned}
\end{align}
for all $\lambda\geq\lambda_0(T+T^2)$.
\end{lm}
Using Lemma \ref{thmm3.1cab}, we  derive the following general Carleman estimate.
\begin{lm}
Let $\mathcal{B}\subset G$ be a nonempty open subset and $d\in\mathbb{R}$. Then, one can find a positive constant $C$ and a large $\lambda_1\geq1$ depending only on $G$, $\mathcal{B}$, $\mu_0$, $\tau_{ij}^0$, $\tau^0$ and $d$ such that for all  $F_1\in L^2_\mathcal{F}(0,T;L^2(G))$, $F\in L^2_\mathcal{F}(0,T;L^2(G;\mathbb{R}^N))$ and $z_T\in L^2_{\mathcal{F}_T}(\Omega;L^2(G))$, the associated solution $(z,Z)$ of \eqref{eqqgbc} satisfies that
\begin{align}\label{3.22carlemgenBack}
\begin{aligned}
\mathcal{I}(d,z)\leq C\bigg[&\lambda^d \mathbb{E} \int_{0}^T \int_{\mathcal{B}} \theta^2 \gamma^d z^2 \, dx \, dt+\lambda^{d-3}\mathbb{E}\iint_Q \theta^2\gamma^{d-3}F_1^2\,dx\,dt\\
&+\lambda^{d-1}\mathbb{E}\iint_Q \theta^2\gamma^{d-1}|F|^2\,dx\,dt+\lambda^{d-1}\mathbb{E}\iint_Q \theta^2\gamma^{d-1}Z^2\,dx\,dt\bigg],
\end{aligned}
\end{align}
for all $\lambda\geq\lambda_1(T+T^2)$.
\end{lm}
\begin{proof}
Set the change of variables
\begin{align}\label{changofvarb}
\widehat{z} = (\lambda\gamma)^{\frac{d-3}{2}} z\quad\textnormal{and}\quad \widehat{Z} = (\lambda\gamma)^{\frac{d-3}{2}} Z.
\end{align}
It is easy to see that
\begin{equation}\label{eqqwih}
		\begin{cases}
			\begin{array}{ll}
			d\widehat{z} +  L_0(t)\widehat{z} \,dt = \left[(\lambda\gamma)^{\frac{d-3}{2}} (F_1+\textnormal{div}(F)) + \frac{d-3}{2} \gamma' \gamma^{-1} \widehat{z}\right] dt + \widehat{Z} \, dW(t)	 & \textnormal{in } Q, \\
				\widehat{z} = 0 & \textnormal{on } \Sigma.
			\end{array}
		\end{cases}
\end{equation}
Using the Carleman estimate \eqref{carback5.8} for the equation \eqref{eqqwih}, we derive that
\begin{align}\label{inn1}
\begin{aligned}
\mathcal{I}(3, \widehat{z}) \leq C \bigg[ &\lambda^3 \mathbb{E} \int_{0}^T \int_{\mathcal{B}} \theta^2 \gamma^3 |\widehat{z}|^2 \, dx \, dt 
+ \mathbb{E} \iint_Q \theta^2 \left| (\lambda\gamma)^{\frac{d-3}{2}} F_1 + \frac{d-3}{2} \gamma' \gamma^{-1} \widehat{z} \right|^2 \, dx \, dt\\
&
+ \lambda^2\mathbb{E}\iint_Q \theta^2\gamma^2 |(\lambda\gamma)^{\frac{d-3}{2}}F|^2 \,dx\,dt+ \lambda^2 \mathbb{E} \iint_Q \theta^2 \gamma^2 |\widehat{Z}|^2 \, dx \, dt \bigg],
\end{aligned}
\end{align}
for any $\lambda\geq C(T+T^2)$. Using that $(a+b)^2\leq 2(a^2+b^2)$ and  $|\gamma' \gamma^{-1}| \leq 3T \gamma$ for the second term on the right-hand side of \eqref{inn1}, we obtain 
\begin{align}\label{eqq2.78}
\begin{aligned}
\mathcal{I}(3, \widehat{z}) \leq C \bigg[ &\lambda^3 \mathbb{E} \int_{0}^T \int_{\mathcal{B}} \theta^2 \gamma^3 |\widehat{z}|^2 \, dx \, dt
+ \lambda^{d-3} \mathbb{E} \iint_Q \theta^2 \gamma^{d-3} F_1^2 \, dx \, dt +T^2\mathbb{E} \iint_Q \theta^2 \gamma^2 |\widehat{z}|^2 \, dx \, dt\\
&+ \lambda^{d-1}\mathbb{E}\iint_Q \theta^2\gamma^{d-1} |F|^2 \,dx\,dt+ \lambda^{2} \mathbb{E} \iint_Q \theta^2 \gamma^{2} |\widehat{Z}|^2 \, dx \, dt \bigg].
\end{aligned}
\end{align}
Recalling \eqref{2.301}, we have 
\begin{align}\label{ineqqafter298}
CT^2\mathbb{E} \iint_Q \theta^2 \gamma^2 |\widehat{z}|^2 \, dx \, dt\leq CT^4\mathbb{E} \iint_Q \theta^2 \gamma^3 |\widehat{z}|^2 \, dx \, dt.
\end{align}
Taking a large $\lambda\geq CT^2$ in \eqref{ineqqafter298}, it follows that
\begin{align*}
CT^2\mathbb{E} \iint_Q \theta^2 \gamma^2 |\widehat{z}|^2 \, dx \, dt\leq \frac{\lambda^2}{2}\mathbb{E} \iint_Q \theta^2 \gamma^3 |\widehat{z}|^2 \, dx \, dt,
\end{align*}
which implies that
\begin{align}\label{ineqqafter298sec}
CT^2\mathbb{E} \iint_Q \theta^2 \gamma^2 |\widehat{z}|^2 \, dx \, dt\leq \frac{\lambda^3}{2}\mathbb{E} \iint_Q \theta^2 \gamma^3 |\widehat{z}|^2 \, dx \, dt.
\end{align}
Combining \eqref{eqq2.78} and \eqref{ineqqafter298sec}, we end up with
\begin{align*}
\begin{aligned}
\mathcal{I}(3, \widehat{z}) \leq C \bigg[ &\lambda^3 \mathbb{E} \int_{0}^T \int_{\mathcal{B}} \theta^2 \gamma^3 |\widehat{z}|^2 \, dx \, dt
+ \lambda^{d-3} \mathbb{E} \iint_Q \theta^2 \gamma^{d-3} F_1^2 \, dx \, dt\\
&+ \lambda^{d-1}\mathbb{E}\iint_Q \theta^2\gamma^{d-1} |F|^2 \,dx\,dt+ \lambda^{2} \mathbb{E} \iint_Q \theta^2 \gamma^{2} |\widehat{Z}|^2 \, dx \, dt \bigg].
\end{aligned}\end{align*}
Finally, recalling \eqref{changofvarb}, we deduce the desired Carleman inequality \eqref{3.22carlemgenBack}.
\end{proof}

\subsection{Main Carleman Inequality}

In this subsection, we derive the following main Carleman estimate for the
coupled system \eqref{adback4.77firstcontro}.

\begin{thm}\label{carlthm32f1.5}
Assume that \eqref{assump1.3} and \eqref{assump1.3B21} hold. Then, there exists
a constant $C>0$, depending only on
$G$, $G_0$, $\widetilde{G}_0$, $\mu_0$, $\tau_{ij}^m$,
$\tau_0$, and $\|a_{21}\|_\infty$, such that for all $h_T,k_T\in L^2_{\mathcal{F}_T}(\Omega;L^2(G))$, the associated solution $(h,k;H,K)$ of the system
\eqref{adback4.77firstcontro} satisfies
\begin{align}\label{carlestcasca1.5sec}
\begin{aligned}
\mathcal{I}(6,h)+\mathcal{I}(3,k)
\leq C\bigg[\lambda^{8}&\mathbb{E}\iint_{Q_0}
\theta^2\gamma^{8}h^2\,dx\,dt+\lambda^{5}\mathbb{E}\iint_Q
\theta^2\gamma^{5}H^2\,dx\,dt\\
&+\lambda^{2}\mathbb{E}\iint_Q
\theta^2\gamma^{2}K^2\,dx\,dt
\bigg],
\end{aligned}
\end{align}
for all sufficiently large
\[
\lambda\geq C\Big[
T+T^2+T^2
\Big(
\|a_{11}\|_\infty^{\frac23}
+\|B_{11}\|_\infty^2
+\|a_{22}\|_\infty^{\frac23}
+\|B_{22}\|_\infty^2
+\|a_{12}\|_\infty^{\frac13}
+\|B_{12}\|_\infty^{\frac12}
\Big)
\Big].
\]
\end{thm}
\begin{proof}
For clarity, we divide the proof into three  steps.\\
\textbf{Step 1. Some notations and preliminaries.}\\
Let \( w_\ell = \theta^2 (\lambda \gamma)^\ell \) with \( \ell \in \mathbb{N} \), and define the subsets \( \widetilde{G}_i \) (\( i = 1, 2 \)) such that
\begin{equation*}
\widetilde{G}_2 \Subset \widetilde{G}_1 \Subset \widetilde{G}_0\subset G_0,
\end{equation*}
where \( \widetilde{G}_0 \) is the set introduced in \eqref{assump1.3}. We also consider the functions \( \rho_i \in C^{\infty}(\mathbb{R}^N) \) so that
\begin{align}\label{assmzeta}
\begin{aligned}
& 0 \leq \rho_i \leq 1, \quad \rho_i = 1 \,\, \text{in} \,\, \widetilde{G}_{3-i}, \quad \text{Supp}(\rho_i) \subset \widetilde{G}_{2-i}, \\ 
& \frac{\Delta \rho_i}{\rho_i^{1/2}} \in L^\infty(G), \quad \frac{\nabla \rho_i}{\rho_i^{1/2}} \in L^\infty(G; \mathbb{R}^N), \quad i=1,2.
\end{aligned}
\end{align}
From \eqref{2.301} and \eqref{assmzeta}, we show that for sufficiently large \( \lambda \geq C(T + T^2) \), one has
\begin{align}\label{esttmforT}
|\partial_t w_\ell|\leq C\lambda^{\ell+2}\theta^2\gamma^{\ell+2},\qquad|\nabla(w_\ell\rho_i)|\leq C\lambda^{\ell+1}\theta^2\gamma^{\ell+1}\rho_i^{1/2},\quad i=1,2.
\end{align}
Applying the Carleman estimate \eqref{3.22carlemgenBack} for the state $h$ of the system \eqref{adback4.77firstcontro}, with \(\mathcal{B} = \widetilde{G}_2\) and \(d = 6\), there exits a constant \(C > 0\) and a large enough \(\lambda\geq C(T+T^2)\) such that
\begin{align}\label{estt4.3fsec1}
    \begin{aligned}
\mathcal{I}(6,h) \leq C \bigg[ &\lambda^6 \mathbb{E} \int_0^T \int_{\widetilde{G}_2} \theta^2 \gamma^6 h^2 \, dx \, dt + \lambda^3 \mathbb{E} \iint_Q \theta^2 \gamma^3 |a_{11}h+a_{21}k|^2 \, dx \, dt \\
& + \lambda^5 \mathbb{E} \iint_Q \theta^2 \gamma^5 |hB_{11}+kB_{21}|^2 \, dx \, dt + \lambda^5 \mathbb{E} \iint_Q \theta^2 \gamma^5 H^2 \, dx \, dt \bigg].
    \end{aligned}
\end{align}
Recalling \eqref{assump1.3B21} and taking a large $\lambda\geq C(T+T^2+\|a_{11}\|_\infty^{\frac{2}{3}}T^2+\|B_{11}\|^2_\infty T^2)$ in \eqref{estt4.3fsec1}, we have that
\begin{align}\label{estt4.3fsec1se}
    \begin{aligned}
\mathcal{I}(6,h) \leq C \bigg[ &\lambda^6 \mathbb{E} \int_0^T \int_{\widetilde{G}_2} \theta^2 \gamma^6 h^2 \, dx \, dt  + \lambda^3 \|a_{21}\|_\infty^2\mathbb{E} \iint_Q \theta^2 \gamma^3 k^2 \, dx \, dt\\
&+ \lambda^5 \mathbb{E} \iint_Q \theta^2 \gamma^5 H^2 \, dx \, dt \bigg].
    \end{aligned}
\end{align}
Using again the Carleman estimate \eqref{3.22carlemgenBack} for the state $k$ of the system \eqref{adback4.77firstcontro}, with \(\mathcal{B} = \widetilde{G}_2\) and \(d = 3\), there exits a constant \(C > 0\) and a large enough \(\lambda\geq C(T+T^2)\) such that
\begin{align}\label{estt4.3fsec}
    \begin{aligned}
\mathcal{I}(3,k) \leq C \bigg[ &\lambda^3 \mathbb{E} \int_0^T \int_{\widetilde{G}_2} \theta^2 \gamma^3 k^2 \, dx \, dt + \mathbb{E} \iint_Q \theta^2 |a_{12}h+a_{22}k|^2 \, dx \, dt \\
& + \lambda^2 \mathbb{E} \iint_Q \theta^2 \gamma^2 |hB_{12}+kB_{22}|^2 \, dx \, dt + \lambda^2 \mathbb{E} \iint_Q \theta^2 \gamma^2 K^2 \, dx \, dt \bigg].
    \end{aligned}
\end{align}
Taking a large $\lambda\geq C(T+T^2+\|a_{22}\|_\infty^{\frac{2}{3}}T^2+\|B_{22}\|^2_\infty T^2)$ in \eqref{estt4.3fsec}, we have that
\begin{align}\label{estt4.3fsec22}
    \begin{aligned}
\mathcal{I}(3,k) \leq C \bigg[ &\lambda^3 \mathbb{E} \int_0^T \int_{\widetilde{G}_2} \theta^2 \gamma^3 k^2 \, dx \, dt + \|a_{12}\|_\infty^2\mathbb{E} \iint_Q \theta^2 h^2 \, dx \, dt\\
& + \lambda^2 \|B_{12}\|^2_\infty\mathbb{E} \iint_Q \theta^2 \gamma^2 h^2 \, dx \, dt + \lambda^2 \mathbb{E} \iint_Q \theta^2 \gamma^2 K^2 \, dx \, dt \bigg].
    \end{aligned}
\end{align}
Combining \eqref{estt4.3fsec1se} and \eqref{estt4.3fsec22}, we get 
\begin{align}\label{ineqqne6.3}
    \begin{aligned}
\mathcal{I}(6,h)+\mathcal{I}(3,k) \leq C \bigg[ &\lambda^6 \mathbb{E} \int_0^T \int_{\widetilde{G}_2} \theta^2 \gamma^6 h^2 \, dx \, dt+\lambda^3 \mathbb{E} \int_0^T \int_{\widetilde{G}_2} \theta^2 \gamma^3 k^2 \, dx \, dt\\
&+ \|a_{12}\|_\infty^2\mathbb{E} \iint_Q \theta^2 h^2 \, dx \, dt + \lambda^3 \|a_{21}\|_\infty^2\mathbb{E} \iint_Q \theta^2 \gamma^3 k^2 \, dx \, dt \\
&+\lambda^2 \|B_{12}\|^2_\infty\mathbb{E} \iint_Q \theta^2 \gamma^2 h^2 \, dx \, dt+ \lambda^5 \mathbb{E} \iint_Q \theta^2 \gamma^5 H^2 \, dx \, dt\\
& + \lambda^2 \mathbb{E} \iint_Q \theta^2 \gamma^2 K^2 \, dx \, dt\bigg],
    \end{aligned}
\end{align}
for large $\lambda\geq C(T+T^2+\|a_{11}\|_\infty^{\frac{2}{3}}T^2+\|B_{11}\|^2_\infty T^2+\|a_{22}\|_\infty^{\frac{2}{3}}T^2+\|B_{22}\|^2_\infty T^2)$. Using again \eqref{estt4.3fsec22}, we obtain
\begin{align}\label{enewestim}
    \begin{aligned}
\mathcal{I}(6,h)+\mathcal{I}(3,k) \leq C \bigg[ &\lambda^6 \mathbb{E} \int_0^T \int_{\widetilde{G}_2} \theta^2 \gamma^6 h^2 \, dx \, dt+\lambda^3 \mathbb{E} \int_0^T \int_{\widetilde{G}_2} \theta^2 \gamma^3 k^2 \, dx \, dt\\
&+ \|a_{12}\|_\infty^2\mathbb{E} \iint_Q \theta^2 h^2 \, dx \, dt +\lambda^2 \|B_{12}\|^2_\infty\mathbb{E} \iint_Q \theta^2 \gamma^2 h^2 \, dx \, dt \\
&+ \lambda^5 \mathbb{E} \iint_Q \theta^2 \gamma^5 H^2 \, dx \, dt+ \lambda^2 \mathbb{E} \iint_Q \theta^2 \gamma^2 K^2 \, dx \, dt\bigg],
    \end{aligned}
\end{align}
for large $\lambda\geq C(T+T^2+\|a_{11}\|_\infty^{\frac{2}{3}}T^2+\|B_{11}\|^2_\infty T^2+\|a_{22}\|_\infty^{\frac{2}{3}}T^2+\|B_{22}\|^2_\infty T^2)$.
To absorb the third and fourth terms on the right-hand side of \eqref{enewestim}, we take $\lambda$ sufficiently large, specifically,
$$\lambda\geq C(T+T^2+\|a_{11}\|_\infty^{\frac{2}{3}}T^2+\|B_{11}\|^2_\infty T^2+\|a_{22}\|_\infty^{\frac{2}{3}}T^2+\|B_{22}\|^2_\infty T^2 +\|a_{12}\|_\infty^{\frac{1}{3}}T^2+\|B_{12}\|^2_\infty T^2
),$$
we then obtain that
\begin{align}\label{enewestimsec}
    \begin{aligned}
\mathcal{I}(6,h)+\mathcal{I}(3,k) \leq C \bigg[ &\lambda^6 \mathbb{E} \int_0^T \int_{\widetilde{G}_2} \theta^2 \gamma^6 h^2 \, dx \, dt+\lambda^3 \mathbb{E} \int_0^T \int_{\widetilde{G}_2} \theta^2 \gamma^3 k^2 \, dx \, dt\\
& + \lambda^5 \mathbb{E} \iint_Q \theta^2 \gamma^5 H^2 \, dx \, dt + \lambda^2 \mathbb{E} \iint_Q \theta^2 \gamma^2 K^2 \, dx \, dt\bigg].
    \end{aligned}
\end{align}
\textbf{Step 2. Absorption of the localized integral associated with the state $k$.}\\
Let us now absorb the second localized integral on the right-hand side of \eqref{enewestimsec}. From the condition \eqref{assump1.3}, and for simplicity,  we assume that
\begin{align*}
a_{21} \geq a_0 > 0 \quad \textnormal{in} \quad (0,T) \times \widetilde{G}_0, \quad \mathbb{P}\textnormal{-a.s.}
\end{align*}
This implies that
\begin{align}\label{firstestforksec}
    a_0 \lambda^{3} \mathbb{E} \int_0^T \int_{\widetilde{G}_2} \theta^2 \gamma^{3} k^2 \, dx \, dt \leq \mathbb{E} \iint_Q w_{3} \rho_1 a_{21} k^2 \, dx \, dt.
\end{align}
Using system~\eqref{adback4.77firstcontro} and applying Itô’s formula to 
\( d(w_{3}\rho_{1}hk) \), we integrate the resulting identity over \(Q\) and 
take the expectation on both sides. Then, by integration by parts and using 
assumption~\eqref{assump1.3B21}, we obtain
\begin{align}\label{forB21}
    \begin{aligned}
    \mathbb{E} \iint_Q w_{3} \rho_1 a_{21} k^2 \, dx \, dt
    &= \mathbb{E} \iint_Q \partial_t w_{3}\, \rho_1 h k \, dx \, dt + \sum_{i,j=1}^N \mathbb{E} \iint_Q \tau_{ij}^1 \frac{\partial (w_{3} \rho_1 k)}{\partial x_i} \frac{\partial h}{\partial x_j} \, dx \, dt \\
    & \quad - \mathbb{E} \iint_Q w_{3} \rho_1 a_{11} h k \, dx \, dt- \mathbb{E} \iint_Q h B_{11}\cdot\nabla(w_{3} \rho_1 k) \, dx \, dt \\
    & \quad + \sum_{i,j=1}^N \mathbb{E} \iint_Q \tau_{ij}^2 \frac{\partial (w_{3} \rho_1 h)}{\partial x_i} \frac{\partial k}{\partial x_j} \, dx \, dt - \mathbb{E} \iint_Q w_{3} \rho_1 a_{12} h^2 \, dx \, dt\\
    & \quad - \mathbb{E} \iint_Q w_{3} \rho_1 a_{22} h k \, dx \, dt  - \mathbb{E} \iint_Q h B_{12}\cdot\nabla(w_{3} \rho_1 h) \, dx \, dt\\
    & \quad - \mathbb{E} \iint_Q k B_{22}\cdot\nabla(w_{3} \rho_1 h) \, dx \, dt.
    \end{aligned}
\end{align}
Let \(\varepsilon>0\) be fixed. By Young's inequality, it is not difficult to see that, for \(\lambda\) sufficiently large, we have
\begin{align}\label{mainestinfirstsec}
\begin{aligned}
 \mathbb{E} \iint_Q w_{3} \rho_1 a_{21} k^2 \, dx \, dt\leq& \,C\varepsilon  \Big[\mathcal{I}(3,k)+\mathcal{I}(6,h)\Big]\\
 &+\frac{C}{\varepsilon}\bigg[\lambda^{7} \mathbb{E} \iint_{Q_0} \theta^2 \gamma^{7}  h^2 \, dx \, dt+\lambda^{5} \mathbb{E} \int_0^T \int_{\widetilde{G}_1} \theta^2 \gamma^{5} |\nabla h|^2 \, dx \, dt\\
 &\hspace{1cm}+\|a_{11}\|_\infty^2\lambda^{3}\mathbb{E} \iint_{Q_0} \theta^2 \gamma^{3}  h^2 \, dx \, dt+\|B_{11}\|_\infty^2\lambda^{5} \mathbb{E} \iint_{Q_0} \theta^2 \gamma^{5}  h^2 \, dx \, dt\\
 &\hspace{1cm}+\|a_{12}\|_\infty \lambda^{3} \mathbb{E} \iint_{Q_0} \theta^2 \gamma^{3}  h^2 \, dx \, dt+\|a_{22}\|_\infty^2\lambda^{3} \mathbb{E} \iint_{Q_0} \theta^2 \gamma^{3}  h^2 \, dx \, dt\\
 &\hspace{1cm}+\|B_{12}\|_\infty^2\lambda^{2} \mathbb{E} \iint_{Q_0} \theta^2 \gamma^{2}  h^2 \, dx \, dt+\|B_{12}\|_\infty\lambda^{4}\mathbb{E} \iint_{Q_0} \theta^2 \gamma^{4} h^2 \, dx \, dt\\
 &\hspace{1cm}+\|B_{22}\|_\infty^2\lambda^{5} \mathbb{E} \iint_{Q_0} \theta^2 \gamma^{5}  h^2 \, dx \, dt +\|B_{22}\|_\infty^2 \lambda^{3} \mathbb{E} \int_0^T \int_{\widetilde{G}_1} \theta^2 \gamma^{3} |\nabla h|^2 \, dx \, dt\\
 &\hspace{1cm}+\lambda^3\mathbb{E} \iint_{Q_0} \theta^2\gamma^3h^2 \,dx\,dt+ \lambda^3\mathbb{E} \iint_{Q_0} \theta^2\gamma^3h^2 \,dx\,dt\bigg].
 \end{aligned}
\end{align} 
By choosing \(\varepsilon>0\) sufficiently small and then taking \(\lambda\) sufficiently large so that
\begin{align}\label{lambdashresec}
\lambda\geq C\Big[T+T^2+T^2\Big(\|a_{11}\|_\infty^{\frac{2}{3}}+\|B_{11}\|^2_\infty +\|a_{22}\|_\infty^{\frac{2}{3}}+\|B_{22}\|^2_\infty +\|a_{12}\|_\infty^{\frac{1}{3}}+\|B_{12}\|_\infty^{\frac{1}{2}}\Big) \Big],
\end{align}
we obtain 
\begin{align}\label{estt4.3fnewonelastsec}
    \begin{aligned}
 \mathcal{I}(6,h) + \mathcal{I}(3,k)\leq C \bigg[& \lambda^{7} \mathbb{E} \iint_{Q_0} \theta^2 \gamma^{7} h^2 \, dx \, dt + \lambda^{5} \mathbb{E} \iint_Q \theta^2 \gamma^{5} H^2 \, dx \, dt \\
&+ \lambda^{2} \mathbb{E} \iint_Q \theta^2 \gamma^{2} K^2 \, dx \, dt + \lambda^{5} \mathbb{E} \int_0^T \int_{\widetilde{G}_1} \theta^2 \gamma^{5} |\nabla h|^2 \, dx \, dt \bigg].
    \end{aligned}
\end{align}
\textbf{Step 3. Absorption of the localized integral associated with $\nabla h$.}\\
We first notice that
\begin{align}\label{esiforgradhsec}
    \tau_0 \lambda^{5} \mathbb{E} \int_0^T \int_{\widetilde{G}_1} \theta^2 \gamma^{5} |\nabla h|^2 \, dx \, dt \leq \mathbb{E} \iint_Q w_{5} \rho_2 \sum_{i,j=1}^N \tau_{ij}^1 \frac{\partial h}{\partial x_i} \frac{\partial h}{\partial x_j} \, dx \, dt.
\end{align}
Computing \( d(w_{5} \rho_2 h^2) \), using the 
assumptions~\eqref{assump1.3B21}, integrating on $Q$ the obtained equality and taking the expectation on both sides, we get that
\begin{align}\label{mainestinfirst2sec}
    \begin{aligned}
    &2 \mathbb{E} \iint_Q w_{5} \rho_2 \sum_{i,j=1}^N \tau_{ij}^1 \frac{\partial h}{\partial x_i} \frac{\partial h}{\partial x_j} \, dx \, dt\\
    &= - \mathbb{E} \iint_Q \partial_t w_{5} \,\rho_2 h^2 \, dx \, dt - 2 \sum_{i,j=1}^N \mathbb{E} \iint_Q \tau_{ij}^1 h \frac{\partial (w_{5} \rho_2)}{\partial x_i} \frac{\partial h}{\partial x_j} \, dx \, dt \\
    & \quad + 2 \mathbb{E} \iint_Q w_{5} \rho_2 a_{11} h^2 \, dx \, dt + 2 \mathbb{E} \iint_Q w_{5} \rho_2 a_{21} h k \, dx \, dt \\
    & \quad + 2\mathbb{E} \iint_Q h B_{11}\cdot\nabla(w_{5} \rho_2 h) \, dx \, dt - \mathbb{E} \iint_Q w_{5} \rho_2 H^2 \, dx \, dt.
    \end{aligned}
\end{align}
For the third term on the right-hand side of \eqref{mainestinfirst2sec}, we apply Young's inequality, we have that
$$2 \mathbb{E} \iint_Q w_{5} \rho_2 a_{11} h^2 \, dx \, dt\leq C\|a_{11}\|_\infty^2\lambda^3\mathbb{E}\iint_Q \theta^2\gamma^3h^2 \,dx\,dt+C\lambda^7\mathbb{E}\iint_Q \theta^2\gamma^7\rho_2^2h^2 \,dx\,dt.$$
Then, taking a large $\lambda\geq C\|a_{11}\|_\infty^{\frac{2}{3}}T^2$, we get
\begin{align}\label{ineqI3estsec}
2 \mathbb{E} \iint_Q w_{5} \rho_2 a_{11} h^2 \, dx \, dt\leq C\lambda^6\mathbb{E}\iint_Q \theta^2\gamma^6\rho_2 h^2 \,dx\,dt+C\lambda^7\mathbb{E}\iint_Q \theta^2\gamma^7\rho_2 h^2 \,dx\,dt.
\end{align}
For the fifth term, we first note that
\begin{align*}
2\mathbb{E} \iint_Q h B_{11}\cdot\nabla(w_{5} \rho_2 h) \, dx \, dt=2\mathbb{E} \iint_Q hB_{11}\cdot(w_5\rho_2\nabla h+h\nabla(w_5\rho_2)) \,dx\,dt.
\end{align*}
It then follows that
\begin{align}\label{ineqqforhb11}
\begin{aligned}
2\mathbb{E} \iint_Q h B_{11}\cdot\nabla(w_{5} \rho_2 h) \, dx \, dt&\leq C\lambda^5\|B_{11}\|_\infty\mathbb{E} \iint_Q \theta^2\gamma^5\rho_2|h||\nabla h| \,dx\,dt\\
&\quad+C\lambda^6\|B_{11}\|_\infty\mathbb{E} \iint_Q \theta^2\gamma^6\rho_2^{1/2}h^2 \,dx\,dt.
\end{aligned}
\end{align} 
For the first term on the right-hand side of \eqref{ineqqforhb11}, by Young's inequality, we have that
\begin{align*}
C\lambda^5\|B_{11}\|_\infty\mathbb{E} \iint_Q \theta^2\gamma^5\rho_2|h||\nabla h| \,dx\,dt\leq &\,C\|B_{11}\|^2_\infty\lambda^5\mathbb{E} \iint_Q \theta^2\gamma^5\rho_2 h^2 \,dx\,dt \\
&+ \tau_0\lambda^5\mathbb{E} \iint_Q \theta^2\gamma^5\rho_2|\nabla h|^2 \,dx\,dt.
\end{align*}
Hence, choosing $\lambda$ sufficiently large so that  $\lambda \ge C \|B_{11}\|_{\infty}^{2} T^{2}$, and recalling \eqref{assmponalpha}, we obtain that
\begin{align}\label{ineqqhgradh}
\begin{aligned}
C\lambda^5\|B_{11}\|_\infty\mathbb{E} \iint_Q \theta^2\gamma^5\rho_2|h||\nabla h| \,dx\,dt\leq&\, C\lambda^6\mathbb{E} \iint_Q \theta^2\gamma^6\rho_2 h^2 \,dx\,dt \\
&+ \mathbb{E} \iint_Q w_{5} \rho_2 \sum_{i,j=1}^N \tau_{ij}^1 \frac{\partial h}{\partial x_i} \frac{\partial h}{\partial x_j} \, dx \, dt.
\end{aligned}
\end{align}
For the second term on the right-hand side of \eqref{ineqqforhb11}, we have
\begin{align*}
C\lambda^6\|B_{11}\|_\infty\mathbb{E} \iint_Q \theta^2\gamma^6\rho_2^{1/2}h^2 \,dx\,dt\leq &\,C\lambda^5\|B_{11}\|^2_\infty\mathbb{E} \iint_{Q} \theta^2 \gamma^{5}\rho_2^{1/2} h^2 \, dx \, dt\\
&+C\lambda^7\mathbb{E} \iint_{Q} \theta^2 \gamma^{7}\rho_2^{1/2} h^2 \, dx \, dt.
\end{align*}
Taking a large $\lambda\geq C\|B_{11}\|^2_\infty T^2$, it follows that
\begin{align}\label{ineqqbef245b11}
\begin{aligned}
C\lambda^6\|B_{11}\|_\infty\mathbb{E} \iint_Q \theta^2\gamma^6\rho_2^{1/2}h^2 \,dx\,dt\leq &\,C\lambda^6\mathbb{E} \iint_{Q} \theta^2 \gamma^{6}\rho_2^{1/2} h^2 \, dx \, dt\\
&+C\lambda^7\mathbb{E} \iint_{Q} \theta^2 \gamma^{7}\rho_2^{1/2} h^2 \, dx \, dt.
\end{aligned}
\end{align}
From \eqref{ineqqforhb11}, \eqref{ineqqhgradh} and \eqref{ineqqbef245b11}, we deduce that
\begin{align}\label{ineqqJ15}
\begin{aligned}
    2\mathbb{E} \iint_Q h B_{11}\cdot\nabla(w_{5} \rho_2 h) \, dx \, dt\leq&\, C\lambda^6\mathbb{E} \iint_{Q_0} \theta^2 \gamma^{6} h^2 \, dx \, dt  + \mathbb{E} \iint_Q w_{5} \rho_2 \sum_{i,j=1}^N \tau_{ij}^1 \frac{\partial h}{\partial x_i} \frac{\partial h}{\partial x_j} \, dx \, dt\\
    &+C\lambda^7\mathbb{E} \iint_{Q_0} \theta^2 \gamma^{7}  h^2 \, dx \, dt.
    \end{aligned}
\end{align}
Let us fix $\varepsilon>0$. Using \eqref{esttmforT} together with \eqref{ineqI3estsec} and \eqref{ineqqJ15}, and applying Young's inequality, it follows that, for $\lambda$ sufficiently large as in \eqref{lambdashresec}, we have
\begin{align}\label{ineqq242}
    \begin{aligned}
\mathbb{E} \iint_Q w_{5} \rho_2 \sum_{i,j=1}^N \tau_{ij}^1 \frac{\partial h}{\partial x_i} \frac{\partial h}{\partial x_j} \, dx \, dt\leq&\,C\varepsilon\Big[\mathcal{I}(6,h)+\mathcal{I}(3,k)\Big]\\
&+\frac{C}{\varepsilon}\lambda^8\mathbb{E} \iint_{Q_0} \theta^2\gamma^8 h^2 \,dx\,dt+C\lambda^6\mathbb{E} \iint_{Q_0} \theta^2\gamma^6 h^2 \,dx\,dt\\
&+C\lambda^7\mathbb{E} \iint_{Q_0} \theta^2\gamma^7 h^2 \,dx\,dt+\lambda^5\mathbb{E} \iint_{Q} \theta^2\gamma^5 H^2 \,dx\,dt\\
&+C\lambda^5\mathbb{E} \iint_{Q_0} \theta^2\gamma^5  h^2 \,dx\,dt.
\end{aligned}
\end{align}
Combining \eqref{estt4.3fnewonelastsec}, \eqref{esiforgradhsec},  \eqref{ineqq242}, and selecting a small \(\varepsilon\) and a large \(\lambda\) such that
$$\lambda\geq C\Big[T+T^2+T^2\Big(\|a_{11}\|_\infty^{\frac{2}{3}}+\|B_{11}\|^2_\infty +\|a_{22}\|_\infty^{\frac{2}{3}}+\|B_{22}\|^2_\infty +\|a_{12}\|_\infty^{\frac{1}{3}}+\|B_{12}\|_\infty^{\frac{1}{2}}\Big) \Big],$$
we end up with
\begin{align*}
\begin{aligned}
\mathcal{I}(6,h) + \mathcal{I}(3,k) \leq C \bigg[ \lambda^{8} \mathbb{E} \iint_{Q_0} \theta^2 \gamma^{8} h^2 \, dx \, dt + \lambda^{5} \mathbb{E} \iint_Q \theta^2 \gamma^{5} H^2 \, dx \, dt+ \lambda^{2} \mathbb{E} \iint_Q \theta^2 \gamma^{2} K^2 \, dx \, dt \bigg].
    \end{aligned}
\end{align*}
This implies the desired Carleman estimate \eqref{carlestcasca1.5sec}, then completing the proof of Theorem \ref{carlthm32f1.5}.
\end{proof}

\section{Observability Results}\label{sec4}

Based on the estimate \eqref{carlestcasca1.5sec}, we first establish the
following observability inequality for \eqref{adback4.77firstcontro}.

\begin{prop}\label{propoibserineforcon}
Assume that \eqref{assump1.3} and \eqref{assump1.3B21} hold. Then, there
exists a constant $C>0$ depending only on $G$, $G_0$, $\widetilde{G}_0$,
$\mu_0$, $\tau_{ij}^m$, $\tau_0$, and $\|a_{21}\|_\infty$ such that for any $h_T, k_T \in L^2_{\mathcal{F}_T}(\Omega;L^2(G))$, the associated solution $(h,k;H,K)$ of the system
\eqref{adback4.77firstcontro} satisfies the following inequality:
\begin{align}\label{obseine4.9cornobin}
\begin{aligned}
\mathbb{E}\int_G \left[|h(0)|^2+|k(0)|^2\right]\,dx
\leq
\exp(CM)\,
\mathbb{E}\iint_Q
\left(h^2\chi_{G_0}+H^2+K^2\right)\,dx\,dt,
\end{aligned}
\end{align}
where the constant $M$ has the following form
\begin{align}\label{observaconst}
    \begin{aligned}
M&=1+T+\frac{1}{T}+\|a_{11}\|_\infty^{\frac{2}{3}}+\|B_{11}\|^2_\infty +\|a_{22}\|_\infty^{\frac{2}{3}}+\|B_{22}\|^2_\infty +\|a_{12}\|_\infty^{\frac{1}{3}}+\|B_{12}\|_\infty^{\frac{1}{2}}\\
&\quad+T(\|a_{11}\|_\infty+\|a_{12}\|_\infty+\|a_{22}\|_\infty+\|B_{11}\|^2_\infty+\|B_{12}\|^2_\infty+\|B_{22}\|^2_\infty).
    \end{aligned}
\end{align}
\end{prop}
\begin{proof}
From the Carleman inequality \eqref{carlestcasca1.5sec}, we deduce that
\begin{align}\label{carles1}
 \begin{aligned}
\lambda^3\mathbb{E}\int_{T/4}^{3T/4}\int_G \theta^2 \gamma^3 (h^2+k^2) \,dx\,dt \leq C \bigg[ &\lambda^8\mathbb{E}\iint_{Q_0} \theta^2 \gamma^{8} h^2 \,dx\,dt + \lambda^8\mathbb{E}\iint_Q \theta^2 \gamma^{8} H^2 \,dx\,dt \\
&+ \lambda^8\mathbb{E}\iint_Q \theta^2 \gamma^{8} K^2 \,dx\,dt \bigg],
\end{aligned}
\end{align}
for large
$$\lambda\geq C\Big[T+T^2+T^2\Big(\|a_{11}\|_\infty^{\frac{2}{3}}+\|B_{11}\|^2_\infty +\|a_{22}\|_\infty^{\frac{2}{3}}+\|B_{22}\|^2_\infty +\|a_{12}\|_\infty^{\frac{1}{3}}+\|B_{12}\|_\infty^{\frac{1}{2}}\Big) \Big].$$
It is not difficult to see that there exists a constant $C>0$ such that for a large $\lambda\geq CT^2$,
$$\lambda^3\theta^2\gamma^3\geq C\exp(-C\lambda T^{-2}),\quad\textnormal{in}\;\;(T/4,3T/4)\times G,$$
$$\lambda^8\theta^2\gamma^8\leq C,\quad\textnormal{in}\;\;Q.$$
Hence, we derive that
\begin{align}\label{carles2f}
 \begin{aligned}
\mathbb{E}\int_{T/4}^{3T/4}\int_G (h^2+k^2) \,dx\,dt \leq C\exp(C\lambda T^{-2}) \mathbb{E} \iint_Q \left( h^2\chi_{G_0} + H^2 + K^2 \right) \, dx \, dt,
\end{aligned}
\end{align}
for any
$$\lambda\geq \lambda_1=C\Big[T+T^2+T^2\Big(\|a_{11}\|_\infty^{\frac{2}{3}}+\|B_{11}\|^2_\infty +\|a_{22}\|_\infty^{\frac{2}{3}}+\|B_{22}\|^2_\infty +\|a_{12}\|_\infty^{\frac{1}{3}}+\|B_{12}\|_\infty^{\frac{1}{2}}\Big) \Big].$$
We now fix $\lambda=\lambda_1$ in \eqref{carles2f}, we get that
\begin{align}\label{carles2}
 \begin{aligned}
\mathbb{E}\int_{T/4}^{3T/4}\int_G (h^2+k^2) \,dx\,dt \leq C\exp(CM_0) \mathbb{E} \iint_Q \left( h^2\chi_{G_0} + H^2 + K^2 \right) \, dx \, dt,
\end{aligned}
\end{align}
where 
$$M_0=1+\frac{1}{T}+\|a_{11}\|_\infty^{\frac{2}{3}}+\|B_{11}\|^2_\infty +\|a_{22}\|_\infty^{\frac{2}{3}}+\|B_{22}\|^2_\infty +\|a_{12}\|_\infty^{\frac{1}{3}}+\|B_{12}\|_\infty^{\frac{1}{2}}.$$
Let \( t \in (0,T) \), by Itô's formula, we have that 
$$ d(h^2 + k^2) = 2h\,dh + 2k\,dk + (dh)^2 + (dk)^2.$$ 
Recalling \eqref{assump1.3B21}, using the equation \eqref{adback4.77firstcontro}, integrating the obtained equality over \( (0,t) \) and taking the expectation on both sides, we get that
\begin{align}\label{ineqq124}
\begin{aligned}
\mathbb{E}\int_0^t\int_G d(h^2+k^2)\,dx
&= 2\mathbb{E}\int_0^t\int_G \sum_{i,j=1}^N \tau_{ij}^1\frac{\partial h}{\partial x_i}\frac{\partial h}{\partial x_j} \,dx\,ds - 2\mathbb{E}\int_0^t\int_G a_{11} h^2 \,dx\,ds  \\
&\hspace{0.4cm} - 2\mathbb{E}\int_0^t\int_G a_{21} h k \,dx\,ds - 2\mathbb{E}\int_0^t\int_G h B_{11}\cdot\nabla h \,dx\,ds \\
&\hspace{0.4cm} + \mathbb{E}\int_0^t\int_G H^2 \,dx\,ds   + 2\mathbb{E}\int_0^t\int_G \sum_{i,j=1}^N \tau_{ij}^2\frac{\partial k}{\partial x_i}\frac{\partial k}{\partial x_j}  \,dx\,ds\\
&\hspace{0.4cm} - 2\mathbb{E}\int_0^t\int_G a_{12} h k \,dx\,ds   - 2\mathbb{E}\int_0^t\int_G a_{22} k^2 \,dx\,ds\\
&\hspace{0.4cm}
- 2\mathbb{E}\int_0^t\int_G h B_{12}\cdot\nabla k \,dx\,ds- 2\mathbb{E}\int_0^t\int_G k B_{22}\cdot\nabla k \,dx\,ds\\
&\hspace{0.4cm}
+ \mathbb{E}\int_0^t\int_G K^2 \,dx\,ds.
\end{aligned}
\end{align}
Using \eqref{assmponalpha} for the first and sixth terms and Young's inequality for the other terms on the right-hand side of \eqref{ineqq124}, we find that for any $\varepsilon>0$,
\begin{align}\label{ineqq124sec}
\begin{aligned}
-\mathbb{E}\int_0^t\int_G d(h^2+k^2)\,dx
&\leq -2\tau_0\mathbb{E}\int_0^t\int_G |\nabla h|^2 \,dx\,ds + 2\|a_{11}\|_\infty\mathbb{E}\int_0^t\int_G  h^2 \,dx\,ds  \\
&\hspace{0.4cm} + 2\|a_{21}\|_\infty\mathbb{E}\int_0^t\int_G  |h| |k| \,dx\,ds + \varepsilon\mathbb{E}\int_0^t\int_G |\nabla h|^2 \,dx\,ds\\
&\hspace{0.4cm}+\frac{C}{\varepsilon}\|B_{11}\|_\infty^2\mathbb{E}\int_0^t\int_G h^2 \,dx\,ds  - \mathbb{E}\int_0^t\int_G H^2 \,dx\,ds \\
&\hspace{0.4cm}  -2\tau_0\mathbb{E}\int_0^t\int_G |\nabla k|^2 \,dx\,ds+ 2\|a_{12}\|_\infty\mathbb{E}\int_0^t\int_G  |h| |k| \,dx\,ds \\
&\hspace{0.4cm} + 2\|a_{22}\|_\infty\mathbb{E}\int_0^t\int_G k^2 \,dx\,ds
+ \varepsilon\mathbb{E}\int_0^t\int_G |\nabla k|^2 \,dx\,ds\\
&\hspace{0.4cm}+\frac{C}{\varepsilon}\|B_{12}\|_\infty^2\mathbb{E}\int_0^t\int_G h^2 \,dx\,ds
+ \varepsilon\mathbb{E}\int_0^t\int_G |\nabla k|^2 \,dx\,ds\\
&\hspace{0.4cm}
+\frac{C}{\varepsilon}\|B_{22}\|_\infty^2\mathbb{E}\int_0^t\int_G k^2 \,dx\,ds- \mathbb{E}\int_0^t\int_G K^2 \,dx\,ds\\
&\hspace{0.4cm}+C\mathbb{E}\int_0^t\int_G h^2 \,dx\,ds+C\mathbb{E}\int_0^t\int_G k^2 \,dx\,ds.
\end{aligned}
\end{align}
Taking a small enough $\varepsilon$, we get that
\begin{align*}
\begin{aligned}
\mathbb{E}\int_G \left[ h^2(0) + k^2(0) \right] \,dx 
&\leq \mathbb{E}\int_G \left[ h^2(t) + k^2(t) \right] \,dx + CM_1\mathbb{E}\int_0^t\int_G \left( h^2 + k^2 \right) \,dx\,ds,
\end{aligned}
\end{align*}
where 
$$M_1=1+\|a_{11}\|_\infty+\|a_{12}\|_\infty+\|a_{22}\|_\infty+\|B_{11}\|^2_\infty+\|B_{12}\|^2_\infty+\|B_{22}\|^2_\infty.$$
Hence, by Gronwall's inequality, it follows that
\begin{align}\label{ineqq237}
\begin{aligned}
\mathbb{E}\int_G \left[ h^2(0) + k^2(0) \right] \,dx 
&\leq \exp(CM_1T)\; \mathbb{E}\int_G \left[ h^2(t) + k^2(t) \right] \,dx.
\end{aligned}
\end{align}
Integrating \eqref{ineqq237} on \( (T/4, 3T/4) \) and combining the obtained result with \eqref{carles2}, we obtain the observability inequality \eqref{obseine4.9cornobin}. This completes the proof of Proposition \ref{propoibserineforcon}.
\end{proof}

From the Carleman estimate \eqref{carlestcasca1.5sec}, it is not difficult to
derive the following unique continuation property for the system
\eqref{adback4.77firstcontro}.

\begin{prop}\label{unicconprop}
Assume that \eqref{assump1.3} and \eqref{assump1.3B21} hold. Then, for any $h_T, k_T \in L^2_{\mathcal{F}_T}(\Omega;L^2(G))$, the associated solution $(h,k;H,K)$ of the system
\eqref{adback4.77firstcontro} satisfies
\begin{align}\label{uniqcontpropercon}
(h,H,K)=(0,0,0)
\quad\textnormal{in } Q_0\times Q^2,
\quad \mathbb{P}\textnormal{-a.s.}
\quad \Longrightarrow \quad
h_T=k_T=0
\quad\textnormal{in } G,
\quad \mathbb{P}\textnormal{-a.s.}
\end{align}
\end{prop}

\section{Proofs of the Main Controllability Results}\label{sec5}
For completeness, this section is dedicated to proving the main results of this paper, as presented in Theorems \ref{thmm1.3insfirst} and \ref{thm1.2approx}. We begin by establishing the null controllability for the system \eqref{ass15cont}.

\begin{proof}[Proof of Theorem \ref{thmm1.3insfirst}]
Let \( y_0, z_0 \in L^2_{\mathcal{F}_0}(\Omega;L^2(G)) \) and \( \xi_1, \xi_2 \in L^2_\mathcal{F}(0,T; L^2(G)) \)  satisfying \eqref{assonxi1xi2}, and define the following linear subspace of the product space 
\( L^2_{\mathcal{F}}(0,T;L^2(G_0)) \times (L^2_{\mathcal{F}}(0,T;L^2(G)))^2 \):
\begin{align*}
\mathcal{Y} = \Big\{ &(\chi_{G_0} h, H, K) \mid (h, k; H, K) \text{ is the solution of } \eqref{adback4.77firstcontro},\text{ for some } h_T, k_T \in L^2_{\mathcal{F}_T}(\Omega; L^2(G)) \Big\}.
\end{align*}
We also consider the following linear functional \(\mathcal{L}\) defined on \(\mathcal{Y}\) as 
\[
\mathcal{L}(\chi_{G_0} h, H, K) = -\mathbb{E} \iint_Q (h \xi_1 + k \xi_2) \, dx \, dt - \mathbb{E}\int_G \left[ h(0) y_0 + k(0) z_0 \right] \, dx.
\]
By Cauchy-Schwarz inequality, we have that
\begin{align*}
|\mathcal{L}(\chi_{G_0} h, H, K)|^2\leq &\,2\bigg[\bigg(\mathbb{E}\iint_Q \theta^{-2}\gamma^{-6}\xi_1^2 \,dx\,dt+\mathbb{E}\iint_Q \theta^{-2}\gamma^{-3}\xi_2^2\,dx\,dt\bigg)\\
&\quad\;\times\bigg(\mathbb{E}\iint_Q \theta^{2}\gamma^{6}h^2 \,dx\,dt+\mathbb{E}\iint_Q \theta^{2}\gamma^{3}k^2\,dx\,dt\bigg)\\
&\quad\;+\bigg(\mathbb{E}\int_G |y_0|^2\,dx+\mathbb{E}\int_G |z_0|^2\,dx\bigg)\\
&\quad\;\times\bigg(\mathbb{E}\int_G |h(0)|^2\,dx+\mathbb{E}\int_G |k(0)|^2\,dx\bigg)\bigg].
\end{align*}
Using the Carleman estimate \eqref{carlestcasca1.5sec} and the observability inequality \eqref{obseine4.9cornobin}, we deduce that 
\begin{align}
    \begin{aligned}
&\,|\mathcal{L}(\chi_{G_0} h, H, K)|^2\\
&\leq\exp(CM)\Big(\|y_0\|^2_{L^2_{\mathcal{F}_0}(\Omega; L^2(G))} + \|z_0\|^2_{L^2_{\mathcal{F}_0}(\Omega; L^2(G))}\\
&\hspace{2.3cm}+\left\|\theta^{-1}\gamma^{-3}\xi_1\right\|^2_{L^2_\mathcal{F}(0,T;L^2(G))}+ \left\|\theta^{-1}\gamma^{-3/2}\xi_2\right\|^2_{L^2_\mathcal{F}(0,T;L^2(G))}\Big)\\
&\qquad\times|(\chi_{G_0} h, H, K)|^2_{L^2_{\mathcal{F}}(0,T;L^2(G_0)) \times (L^2_{\mathcal{F}}(0,T;L^2(G)))^2},
    \end{aligned}
\end{align}
where $M$ is the constant defined in \eqref{observaconst}. Therefore, the functional \(\mathcal{L}\) is bounded on \(\mathcal{Y}\), and we have that
\begin{align}\label{esthanbancont}
\begin{aligned}
|\mathcal{L}|_{\mathcal{L}(\mathcal{Y}; \mathbb{R})} 
\leq \sqrt{\exp(CM)}\Big(&\|y_0\|_{L^2_{\mathcal{F}_0}(\Omega; L^2(G))} + \|z_0\|_{L^2_{\mathcal{F}_0}(\Omega; L^2(G))}\\
&+\left\|\theta^{-1}\gamma^{-3}\xi_1\right\|_{L^2_\mathcal{F}(0,T;L^2(G))}+ \left\|\theta^{-1}\gamma^{-3/2}\xi_2\right\|_{L^2_\mathcal{F}(0,T;L^2(G))}\Big).
\end{aligned}
\end{align}
By the Hahn–Banach theorem, we can extend \(\mathcal{L}\) to a bounded linear functional on the whole space 
\( L^2_{\mathcal{F}}(0,T; L^2(G_0)) \times (L^2_{\mathcal{F}}(0,T; L^2(G)))^2\). For simplicity, we still denote this extension by \(\mathcal{L}\). 
Using  the Riesz representation theorem, there exist controls 
\[ (f, g_1, g_2) \in L^2_{\mathcal{F}}(0,T; L^2(G_0)) \times (L^2_{\mathcal{F}}(0,T; L^2(G)))^2, \]
such that
\begin{align}\label{fireq11cont}
\begin{aligned}
&\,-\mathbb{E} \iint_Q (h \xi_1 + k \xi_2) \, dx \, dt - \mathbb{E}\int_G \left[ h(0) y_0 + k(0) z_0 \right] \, dx 
\\
&= \mathbb{E} \iint_Q (f \chi_{G_0}  h + g_1 H + g_2 K) \, dx \, dt.
\end{aligned}
\end{align}
Moreover, by \eqref{esthanbancont}, we have 
\begin{align*}
    \begin{aligned}
    &\, \|f\|_{L^2_\mathcal{F}(0,T; L^2(G_0))} + \|g_1\|_{L^2_\mathcal{F}(0,T; L^2(G))} + \|g_2\|_{L^2_\mathcal{F}(0,T; L^2(G))} \\
    &\leq \sqrt{\exp(CM)} \bigg( \|y_0\|_{L^2_{\mathcal{F}_0}(\Omega; L^2(G))} + \|z_0\|_{L^2_{\mathcal{F}_0}(\Omega; L^2(G))} \\
&\hspace{3cm}+\left\|\theta^{-1}\gamma^{-3}\xi_1\right\|_{L^2_\mathcal{F}(0,T;L^2(G))} + \left\|\theta^{-1}\gamma^{-3/2}\xi_2\right\|_{L^2_\mathcal{F}(0,T;L^2(G))} \bigg).
    \end{aligned}
\end{align*}
We claim that the obtained controls \( f \), \( g_1 \), and \( g_2 \) are the desired controls. Indeed, for any \( h_T, k_T \in L^2_{\mathcal{F}_T}(\Omega; L^2(G)) \), by applying Itô's formula to the systems \eqref{ass15cont} and \eqref{adback4.77firstcontro}, we notice that
\begin{align}\label{eqq4.4it}
d(yh + zk) = y \, dh + h \, dy + z \, dk + k \, dz + dy\,dh + dz\,dk.
\end{align}
Integrating \eqref{eqq4.4it} over $Q$ and taking the expectation on both sides, and then applying Fubini’s theorem, we obtain
\begin{align}\label{eqq1.2cor}
\begin{aligned}
\mathbb{E} \int_G \left[ y(T) h_T + z(T) k_T \right] \, dx
= & \, \mathbb{E} \iint_Q \left( h \xi_1 + k \xi_2 + f \chi_{G_0}  h + g_1 H + g_2 K \right) \, dx \, dt \\
& + \mathbb{E}\int_G \left[ h(0) y_0 + k(0) z_0 \right] \, dx.
\end{aligned}
\end{align}
Combining \eqref{eqq1.2cor} and \eqref{fireq11cont}, we deduce that
\[
\mathbb{E} \int_G \left[ y(T) h_T + z(T) k_T \right] \, dx = 0.
\]
Since \( h_T \) and \( k_T \) can be chosen arbitrarily in \( L^2_{\mathcal{F}_T}(\Omega; L^2(G)) \), it follows that
\[
y(T,\cdot) = z(T,\cdot) = 0 \quad \text{in } G, \quad \mathbb{P}\text{-a.s.}
\]
This concludes the proof of Theorem \ref{thmm1.3insfirst}.
\end{proof}

Let us now prove the approximate controllability result for the system \eqref{ass15cont}.
\begin{proof}[Proof of Theorem \ref{thm1.2approx}]
Let \( \xi_1, \xi_2 \in L^2_\mathcal{F}(0,T; L^2(G)) \) and \( y_0, z_0 \in L^2_{\mathcal{F}_0}(\Omega;L^2(G)) \). To establish the approximate controllability of \eqref{ass15cont}, it is equivalent to prove that the following set of reachable states
\begin{align*}
\mathcal{R}_T=\Big\{&(y(T),z(T)) \,\,|\,\, (y,z) \;\;\textnormal{is solution to }\; \eqref{ass15cont}\;\;\textnormal{with some controls}\\
&\,(f,g_1,g_2)\in L^2_\mathcal{F}(0,T; L^2(G_0)) \times (L^2_\mathcal{F}(0,T; L^2(G)))^2 \Big\}
\end{align*}
is dense in \( L^2_{\mathcal{F}_T}(\Omega; L^2(G))\times L^2_{\mathcal{F}_T}(\Omega; L^2(G)) \). By contradiction argument, we suppose that there exist a nonzero element \( (\eta_1, \eta_2) \in L^2_{\mathcal{F}_T}(\Omega; L^2(G))\times L^2_{\mathcal{F}_T}(\Omega; L^2(G)) \) such that
\begin{equation}\label{equa5.5app}
    \mathbb{E} \int_G \left[ y(T) \eta_1 + z(T) \eta_2 \right] \, dx = 0 \quad \text{for all }\;\; (y(T), z(T)) \in \mathcal{R}_T.
\end{equation}
Similarly to \eqref{eqq1.2cor}, applying Itô's formula to the systems \eqref{ass15cont} and \eqref{adback4.77firstcontro} (with \(h_T \equiv \eta_1\) and \(k_T \equiv \eta_2\)) yields
\begin{align}\label{eqq1.2corapp}
\begin{aligned}
\mathbb{E} \int_G \left[ y(T) \eta_1 + z(T) \eta_2 \right] \, dx 
= & \, \mathbb{E} \iint_Q \left( h \xi_1 + k \xi_2 + f \chi_{G_0}  h + g_1 H + g_2 K \right) \, dx \, dt \\
& + \mathbb{E}\int_G \left[ h(0) y_0 + k(0) z_0 \right] \, dx.
\end{aligned}
\end{align}
Combining \eqref{eqq1.2corapp} and \eqref{equa5.5app}, we conclude that
\begin{align*}
\mathbb{E} \iint_Q \left( h \xi_1 + k \xi_2 + f \chi_{G_0}  h + g_1 H + g_2 K \right) \, dx \, dt + \mathbb{E}\int_G \left[ h(0) y_0 + k(0) z_0 \right] \, dx = 0.
\end{align*}
Choosing \( \xi_1 \equiv \xi_2 \equiv y_0 \equiv z_0\equiv 0 \), we obtain that
\begin{equation*}
\mathbb{E} \iint_Q \left( f \chi_{G_0}  h + g_1 H + g_2 K \right) \, dx \, dt = 0,
\end{equation*}
which holds for all \( (f, g_1, g_2) \in L^2_\mathcal{F}(0,T; L^2(G_0)) \times (L^2_\mathcal{F}(0,T; L^2(G)))^2 \). Therefore, we conclude that
\[
(h, H, K) = (0, 0, 0) \quad \text{in} \; Q_0 \times Q^2, \quad \mathbb{P}\textnormal{-a.s.}
\]
By Proposition \ref{unicconprop}, this implies that \( \eta_1 \equiv \eta_2 \equiv 0 \) in \( G \), \( \mathbb{P}\textnormal{-a.s.} \), which is a contradiction. Therefore, the set \( \mathcal{R}_T \) is dense in \( L^2_{\mathcal{F}_T}(\Omega; L^2(G))\times L^2_{\mathcal{F}_T}(\Omega; L^2(G)) \), completing the proof of Theorem \ref{thm1.2approx}.
\end{proof}

\begin{center}
\end{center}
\begin{center}
\end{center}

\begin{flushleft}
$^\dagger$ Dipartimento di Matematica, Università degli Studi di Salerno,  
Via Giovanni Paolo II, 132, 84084 Fisciano (SA), Italy.  
\\\quad Email: \texttt{aelgrou@unisa.it}, \texttt{fgregorio@unisa.it}, \texttt{arhandi@unisa.it}  \\
$^*$Corresponding author.
\end{flushleft}

\begin{thebibliography}{10}

\bibitem{allafrasl22}
B. Allal, G. Fragnelli, and J. Salhi. 
\newblock Null controllability for degenerate parabolic equations with a nonlocal space term.
\newblock{\em  Discrete and Continuous Dynamical Systems - S}, \textbf{17} (2024),  1821--1856. 


\bibitem{surveyAmmarKBGT}
F. Ammar-Khodja, A. Benabdallah, M. González Burgos, and L. de Teresa.  
\newblock Recent results on the controllability of linear coupled parabolic problems: a survey.
\newblock{\em Mathematical Control and Related Fields}, \textbf{1} (2011), 267--306.

\bibitem{elgconvec23}
M. Baroun, S. Boulite, A. Elgrou, and L. Maniar.
\newblock Null controllability for stochastic parabolic equations coupled by first and zero order terms. 
\newblock {\em Applied Mathematics \& Optimization}, \textbf{91} (2025), 31.



\bibitem{Preprielgr24}
S. Boulite, A. Elgrou, and L. Maniar.
\newblock Null controllability for cascade systems of coupled backward stochastic parabolic equations with one distributed control. 
\newblock {\em Journal of Mathematical Analysis and Applications}, \textbf{549} (2025), 129489.


\bibitem{withouextra}
S. Boulite, A. Elgrou, L. Maniar, and O. Oukdach.  
\newblock Controllability for forward stochastic parabolic equations with dynamic boundary conditions without extra forces.
\newblock {\em International Journal of Control}, \textbf{98} (2024), 1913–1923.




\bibitem{Carl39}
T. Carleman. 
\newblock Sur un problème d’unicité pour les systèmes d’équations aux dérivées partielles à deux variables indépendantes.
\newblock {\em Arkiv för Matematik, Astronomi och Fysik}, \textbf{26} (1939), 1--9.


\bibitem{coron07}
J.-M. Coron.
\newblock{\em Control and nonlinearity}. 
\newblock{ American Mathematical Society}, 2007.




\bibitem{evans}
L. C. Evans.
\newblock{\em Partial differential equations}. 
\newblock{ American Mathematical Society}, 2022.



\bibitem{Fadilicasc}
M. Fadili.
\newblock Controllability of a backward stochastic cascade system of coupled parabolic heat equations by one control force. 
\newblock{\em Evolution Equations and Control Theory}, \textbf{12} (2023), 446--458.

\bibitem{ferluzua16}
E. Fern{\'a}ndez-Cara, Q. L{\"u} and E. Zuazua. 
\newblock Null controllability of linear heat and wave equations with nonlocal spatial
terms. 
\newblock{\em SIAM Journal on Control and Optimization}, \textbf{54} (2016), 2009–2019.

\bibitem{fernandez2006global}
E. Fern{\'a}ndez-Cara and S. Guerrero.
\newblock Global Carleman inequalities for parabolic systems and applications to controllability.
\newblock {\em SIAM Journal on Control and Optimization}, \textbf{45} (2006), 1395--1446.

 



 
\bibitem{fursikov1996controllability}
A. V. Fursikov and O. Yu. Imanuvilov.
\newblock Controllability of evolution equations. 
\newblock{\em Lecture Note Series 34, Research Institute of Mathematics, Seoul National University,} (1996).





\bibitem{GonBurTer}
M. González-Burgos and L. de Teresa. 
\newblock Controllability results for cascade systems of $m$-coupled parabolic PDEs by one control force. 
\newblock{\em Portugaliae Mathematica}, \textbf{67} (2010), 91--113.




\bibitem{gureSiam07}
S. Guerrero.
\newblock Null controllability of some systems of two parabolic equations with one control force. 
\newblock {\em SIAM Journal on Control and Optimization}, \textbf{46} (2007), 379--394.


\bibitem{hupeng91}
Y. Hu and S. Peng.
\newblock Adapted solution of a backward semilinear stochastic evolution equation. 
\newblock {\em Stochastic Analysis and Applications}, \textbf{9} (1991), 445--459.





\bibitem{LiQi}	
H. Li and Q. L{\"u}. 
\newblock Null controllability for some systems of two backward stochastic heat equations with one control force. 
\newblock{\em Chinese Annals of Mathematics, Series B}, \textbf{33} (2012), 909--918.


\bibitem{lions1972some}
J.-L. Lions.
\newblock Optimal control of systems governed by partial differential equations. 
\newblock{\em Springer-Verlag, New York-Berlin}, (1971).





    
	
	
\bibitem{liu2014global}
X. Liu.
\newblock Global Carleman estimate for stochastic parabolic equations, and its application.
\newblock{\em ESAIM: Control, Optimisation and Calculus of Variations}, \textbf{20} (2014), 823--839.

 

\bibitem{liu14couplfor}
X. Liu.
\newblock Controllability of some coupled stochastic parabolic systems with fractional order spatial differential operators by one control in the drift. 
\newblock{\em SIAM Journal on Control and Optimization}, \textbf{52} (2014), 836--860.

\bibitem{LiuuLiuX}
 L. Liu and X. Liu.
 \newblock Controllability and observability of some coupled stochastic parabolic systems. 
 \newblock{\em Mathematical Control and Related Fields}, \textbf{8} (2018), 829--854.
 

\bibitem{lotorovozski}
S. V. Lototsky and  B. L. Rozovsky.
\newblock{\em  Stochastic partial differential equations}. Cham: Springer, 2017.



\bibitem{lu2011some}
Q. L{\"u}.
\newblock Some results on the controllability of forward stochastic heat equations with control on the drift.
\newblock{\em Journal of Functional Analysis}, \textbf{260} (2011), 832--851.

\bibitem{luZhang22mcrf}
Q. L{\"u} and X. Zhang.
\newblock A concise introduction to control theory for stochastic partial differential equations.
\newblock{\em Mathematical Control and Related Fields}, \textbf{12} (2022), 847--954.



\bibitem{lu2021mathematical}
Q. L{\"u} and X. Zhang.
\newblock{\em Mathematical control theory for stochastic partial differential equations. } 
\newblock{ Probability Theory and Stochastic Modelling, 101. Springer, Cham}, 2021.



\bibitem{obeman255}
O. Oukdach, S. Boulite, A. Elgrou, and L. Maniar. 
\newblock Stackelberg–Nash null controllability for stochastic parabolic equations. 
\newblock{\em Mathematical Methods in the Applied Sciences}, \textbf{48} (2025), 13164--13176.



\bibitem{tang2009null}
S. Tang and X. Zhang.
\newblock Null controllability for forward and backward stochastic parabolic equations.
\newblock {\em SIAM Journal on Control and Optimization}, \textbf{48} (2009), 2191--2216.



 







\bibitem{Yamam2009invePrb}
M. Yamamoto. 
\newblock Carleman estimates for parabolic equations and applications.
\newblock{\em Inverse Problems}, \textbf{25} (2009), 123013.





\bibitem{yansun2011}
Y. Yan and F. Sun.
\newblock Insensitizing controls for a forward stochastic heat equation. 
\newblock{\em  Journal of Mathematical Analysis and Applications}, \textbf{384} (2011), 138--150.



 



\bibitem{observineqback}
D. Yang and J. Zhong. 
\newblock Observability inequality of backward stochastic heat equations for measurable sets and its applications. 
\newblock{\em SIAM Journal on Control and Optimization}, \textbf{54} (2016), 1157--1175.


\bibitem{wangNull24}
Y. Wang.
\newblock Null controllability for stochastic coupled systems of fourth order parabolic equations.
\newblock{\em Journal of Mathematical Analysis and Applications}, \textbf{538} (2024), 128426.

\bibitem{Zab95}
J. Zabczyk. 
\newblock{\em Mathematical control theory. }
\newblock{ Springer International Publishing,} 2020.




\bibitem{zuazua2007controllability}
E. Zuazua.
\newblock Controllability and observability of partial differential equations: some results and open problems.
\newblock{\em Handbook of Differential Equations: Evolutionary Equations}, \textbf{3} (2007), 527--621.







 \end{thebibliography}
\end{document}